\documentclass[psamsfonts,reqno]{amsart}
\usepackage{amssymb,eucal,latexsym,xy,graphics}
\xyoption{all}
\xyoption{ps}

\theoremstyle{plain}

\newtheorem{lemma}{Lemma}[section]
\newtheorem{theorem}[lemma]{Theorem}
\newtheorem{proposition}[lemma]{Proposition}
\newtheorem{corollary}[lemma]{Corollary}
\newtheorem{example}[lemma]{Example}
\newtheorem{claim}{Claim}

\newtheorem*{stat}{\name}
\newcommand{\name}{testing}

\theoremstyle{definition}
\newtheorem{definition}[lemma]{Definition}
\newtheorem{problem}{Problem}

\theoremstyle{remark}
\newtheorem{remark}[lemma]{Remark}

\newcommand{\qedc}{{\qed}~{\rm Claim~{\theclaim}.}}

\newenvironment{cproof}{\begin{proof}[Proof of Claim.]}
{\qedc\renewcommand{\qed}{}\end{proof}}

\newcommand{\case}[1]%
{\smallskip\noindent\textbf{\textit{Case}\ {#1}.}}

\numberwithin{equation}{section}
\numberwithin{figure}{section}

\newcommand{\pup}[1]{\textup{(}{#1}\textup{)}}

\DeclareMathOperator{\Bool}{Bool}
\DeclareMathOperator{\Retr}{Retr}
\DeclareMathOperator{\Ob}{Ob}
\DeclareMathOperator{\dom}{dom}
\DeclareMathOperator{\im}{im}

\DeclareMathOperator{\Con}{Con}
\DeclareMathOperator{\Id}{Id}

\DeclareMathOperator{\Conc}{Con_c}
\DeclareMathOperator{\lh}{lh}
\DeclareMathOperator{\J}{J}
\DeclareMathOperator{\At}{At}
\DeclareMathOperator{\NAt}{NAt}
\DeclareMathOperator{\M}{M}

\newcommand{\iso}{\mathrm{iso}}

\newcommand{\tr}{\vartriangleleft}

\newcommand{\utr}{\trianglelefteq}
\newcommand{\nutr}{\ntrianglelefteq}
\newcommand{\res}{\mathbin{\restriction}}

\newcommand{\cA}{\mathcal{A}}
\newcommand{\cB}{\mathcal{B}}
\newcommand{\cC}{\mathcal{C}}
\newcommand{\cD}{\mathcal{D}}

\newcommand{\cR}{\mathcal{R}}
\newcommand{\cS}{\mathcal{S}}
\newcommand{\cM}{\mathcal{M}}
\newcommand{\cI}{\mathcal{I}}
\newcommand{\cJ}{\mathcal{J}}
\newcommand{\cL}{\mathcal{L}}
\newcommand{\ZZ}{\mathbb{Z}}
\newcommand{\Rep}{\mathcal{R}_{\mathrm{ep}}}
\newcommand{\AM}[1]{\cA\res\nobreak{#1}}
\newcommand{\AMs}[1]{\cA\res^*\nobreak{#1}}
\newcommand{\AMt}[1]{\widetilde{\cA}\res\nobreak{#1}}
\newcommand{\Pow}{\mathfrak{P}}
\newcommand{\at}{\mathrm{at}}

\newcommand{\jirr}{join-ir\-re\-duc\-i\-ble}
\newcommand{\mirr}{meet-ir\-re\-duc\-i\-ble}

\newcommand{\jsd}{join-sem\-i\-dis\-trib\-u\-tive}

\newcommand{\eps}{\varepsilon}
\newcommand{\epst}{\tilde{\varepsilon}}
\newcommand{\es}{\varnothing}
\newcommand{\into}{\hookrightarrow}
\newcommand{\onto}{\twoheadrightarrow}

\newcommand{\seq}[1]{\left\langle{#1}\right\rangle}
\newcommand{\seqm}[2]{\seq{#1\mid#2}}
\newcommand{\famm}[2]{\left(#1\mid#2\right)}
\newcommand{\set}[1]{\{#1\}}
\newcommand{\setm}[2]{\set{#1\mid#2}}

\newcommand{\go}{\omega}

\newcommand{\xB}{\mathbf{B}}
\newcommand{\Ti}{\mathbf{Ti}}
\newcommand{\GS}{\mathbf{GS}}
\newcommand{\ol}[1]{\overline{#1}}

\newcommand{\two}{\mathbf{2}}
\newcommand{\three}{\mathbf{3}}
\newcommand{\dnw}{\mathbin{\downarrow}}

\newcommand{\id}{\mathrm{id}}
\newcommand{\jz}{$\langle\vee,0\rangle$}
\newcommand{\jzu}{$\langle\vee,0,1\rangle$}
\newcommand{\jzs}{\jz-semi\-lat\-tice}
\newcommand{\jzus}{\jzu-semi\-lat\-tice}
\newcommand{\jzh}{\jz-ho\-mo\-mor\-phism}

\newcommand{\jzuh}{\jzu-ho\-mo\-mor\-phism}
\newcommand{\jze}{\jz-em\-bed\-ding}
\newcommand{\jzue}{\jzu-em\-bed\-ding}

\begin{document}

\title[Functorial inverses without adjunction]%
{Distributive semilattices as retracts of ultraboolean ones;
functorial inverses\\ without adjunction}

\author[F.~Wehrung]{Friedrich Wehrung}
\address{LMNO, CNRS UMR 6139\\ D\'epartement de Math\'ematiques\\
Universit\'e de Caen\\ 14032 Caen Cedex\\ France}
\email{wehrung@math.unicaen.fr}

\urladdr{http://www.math.unicaen.fr/\~{}wehrung}

\date{\today}



\thanks{The author was partially supported by the institutional
grant CEZ:J13/98:1132000007a, by FRVS 2731/2003, and by the Fund
of Mobility of the Charles University (Prague)}

\dedicatory{Souviens toi, ma petite Lynn,\\
la perle du petit dragon\dots}

\begin{abstract}
A \jzs\ is \emph{ultraboolean}, if it is a directed union of finite
Boolean \jzs s. We prove that every distributive \jzs\ is a
retract of some ultraboolean \jzs. This is established by proving
that every finite distributive \jzs\ is a retract of some finite
Boolean \jzs, and this in a \emph{functorial} way. This result is,
in turn, obtained as a particular case of a category-theoretical
result that gives sufficient conditions, for a functor $\Pi$, to
admit a right inverse. The particular functor $\Pi$ used for the
abovementioned result about ultraboolean semilattices has neither a
right nor a left adjoint.
\end{abstract}	

\maketitle

\section{Introduction}\label{S:Intro}

Our general kind of problem is the following. We are given a
functor $\mathbf{F}$ from a category $\cA$ to a
category $\cB$, we wish to investigate
whether $\mathbf{F}$ has a right inverse (up to equivalence).
Also, we suppose that we know how to do this on a subcategory
of~$\cB$, or, more generally, on a given class of diagrams
of~$\cB$. We wish to set a general framework that will enable us,
under certain conditions, to find a right inverse of~$\mathbf{F}$
on a much larger class of diagrams of~$\cB$.

How to do this will be stated precisely in a further
paper~\cite{RetLift}. The present paper is intended to provide a
start for that program, and it is motivated by the following
example. We denote by $\cL$ the category of all lattices, by $\cD$
the category of all distributive \jzs s, and by
$\Conc\colon\cL\to\cD$ the functor that with a lattice~$L$
associates its semilattice $\Conc L$ of compact congruences,
extended naturally to lattice homomorphisms. It is
a well-known open problem, stated by R.\,P. Dilworth in 1945,
whether every distributive \jzs\ is isomorphic to $\Conc L$ for
some lattice $L$. We wish to reduce that problem, or
rather some stronger versions about diagrams of
semilattices, to a smaller class of distributive \jzs s for which
calculations are easier. Our candidate is the following.

\begin{definition}\label{D:Ultrabool}
A \jzs\ is \emph{ultraboolean}, if
it is a directed union of finite Boolean \jzs s.
\end{definition}

Hence every ultraboolean \jzs\ is distributive (the converse is
trivially false, see the three-element chain).

For the present paper's needs, everything boils down to expressing
members of the larger class (distributive semilattices) as
\emph{retracts} of the members of the smaller class (ultraboolean
semilattices). Furthermore, such a retraction needs to be
\emph{functorial}. We shall refer to this problem as the
\emph{ultraboolean retraction problem}. At first sight, it is not
clear whether the functoriality restriction might cause a problem.
Indeed, every finite distributive lattice~$D$ is a retract of a
finite Boolean lattice~$B$. For example, as in
\cite[Section~1]{Ruzi}, we can embed $D$ into the power set
$B=\Pow(\J(D))$, where $\J(D)$ denotes the poset (i.e., partially
ordered set) of \jirr\ elements of $D$, \emph{via} the map
 \[
 a\mapsto\setm{p\in\J(D)}{p\leq a}.
 \]
This map has a retraction, given by $X\mapsto\bigvee X$. As in
\cite[Section~1]{Ruzi}, one can extend `canonically' any \jze\
$f\colon D\into E$ to a \jzh\ $g\colon\Pow(\J(D))\to\Pow(\J(E))$;
however, even for $f=\id_D$, the map $g$ might not be an embedding!
Hence this `functor' preserves neither monomorphisms nor, in fact,
identities, and thus it is not sufficient to solve the ultraboolean
retraction problem.

In order to solve that problem, we need to embed any finite
distributive \jzs\ $D$ into some finite Boolean \jzs\ $\Phi(D)$,
\emph{via} a \jze\ $\eps_D\colon D\into\Phi(D)$, with a retraction
$\mu_D\colon\Phi(D)\onto D$, these data being subjected to
functoriality conditions, stated precisely in
Section~\ref{S:FuncRet}. Here are some caveats:

\begin{itemize}
\item[---] Solving the problem `without the retraction' $\mu_D$
is easy: namely, embed~$D$ into the universal Boolean semilattice
$\Bool(D)$ over $D$. For this construction, the
corresponding embedding $\eps_D\colon D\into\Bool(D)$ is
\emph{not} a meet-embedding as a rule. An explicit construction
is given by $\Bool(D)=\Pow(D^=)$ (where $D^==D\setminus\set{1}$),
$\eps_D(a)=\setm{x\in D^=}{a\nleq x}$ (for all $a\in D$).
On the other hand, any
\jze\ $f\colon D\into E$ is turned to a \emph{lattice} embedding
$g\colon\Bool(D)\into\Bool(E)$! However, the retracts are lost,
for the canonical retraction from $\Bool(D)$ onto~$D$ does not
satisfy the required commutation conditions.

\item[---] For a finite distributive \jzs\ $D$,
the canonical map from $D$ into $\Pow(\J(D))$ is, in fact, a
\emph{lattice} embedding. However, the requirement that all the
maps $\eps_D\colon D\into\Phi(D)$ be lattice embeddings is too
strong to solve the ultraboolean retraction problem. This is
showed by a counterexample in Section~\ref{S:SimultEmb}.
\end{itemize}

Nevertheless, we prove that the ultraboolean retraction problem
has a positive solution. This result is, actually, an immediate
application of a more general categorical principle, stated in
Theorem~\ref{T:Main}. This principle states sufficient conditions
for every object of a category~$\cA$ to be a retract of some
object of a category~$\cB$, and this functorially. Although some
aspects of the formulations might remind of the Adjoint Functor
Theorem, it is not hard to prove that in the particular case of
the ultraboolean retraction problem, the functorial inverse that
we construct does not arise from a functorial adjunction, see
Proposition~\ref{P:NoAdj}.

The importance of finite, simple, atomistic lattices for
representation problems is highlighted in the paper of P.\,P.
P\'alfy and P. Pudl\'ak \cite{P5}, where it is proved that if a
finite, simple lattice $L$ whose atoms join to the unit is
isomorphic to the congruence lattice of a finite algebra, then it
isomorphic to the congruence lattice of a finite set with a finite
group action. With this in mind and by using a trick of G.
Gr\"atzer and E.\,T. Schmidt, we give, in
Section~\ref{S:Tisch}, an easy proof of the result that every
\jzs\ is a retract of some directed \jz-union of finite,
(lattice-)simple, atomistic lattices, and this in a functorial
way. Although this proof does not use the result of
Theorem~\ref{T:Main}, further potential uses of
Theorem~\ref{T:Main} are suggested by open problems such as
Problem~\ref{Pb:SD+} (see Section~\ref{S:Pbs}).

While the present paper deals with the \emph{existence} of
functorial retractions, the paper~\cite{RetLift} deals with how to
\emph{use} functorial retractions in order to prove that certain
functors have large range.

While this paper is mainly category-theoretical, it aims at
building up tools that will be used later in universal algebra. For
this reason, the author chose to write it in probably more detail
than a category theorist would wish, with the hope to make it
reasonably intelligible to members of both communities.

However, a direct semilattice-theoretical proof of
Theorem~\ref{T:FRSemil} (solution of the ultraboolean retraction
problem) is not easier than the categorical proof
involving Theorem~\ref{T:Main}, and it does not lead itself to
further potential generalizations such as those suggested in
Section~\ref{S:Pbs}. This, together with the categorical approach
required in \cite{RetLift}, motivates our choice of the
language of categories instead of the one of universal algebra.

\section{Basic concepts}\label{S:Basic}

Most of our categorical notions are borrowed from S. Mac
Lane~\cite{McLa}. For a category $\cC$, we shall denote by $\Ob\cC$
the class of \emph{objects} of $\cC$, by $\cC^{\iso}$ the category
whose objects are those of $\cC$ and whose morphisms are the
isomorphisms of~$\cC$. We shall denote by $\dom f$ the domain of a
morphism $f$ of $\cC$. As usual, a morphism in~$\cC$ is a
\emph{monic} (resp., a \emph{section}), if it is left cancellable
(resp., left invertible) for the composition of morphisms. Of
course, every section is a monic.

We shall view every quasi-ordered set $\seq{P,\utr}$ as a category
in which hom-sets have at most one element. Technically speaking,
our quasi-ordered sets may be proper classes, but in our context
this will create no difficulty. For $p\utr q$ in $P$, we shall
denote by $p\to q$ the unique morphism from $p$ to $q$. An
\emph{ideal} of $P$ is a subset~$X$ of $P$ such that $p\utr x$
implies that $p\in X$, for all $\seq{p,x}\in P\times X$. We denote
by~$\dnw X$ the ideal generated by $X$, for all
$X\subseteq P$, and we put $\dnw p=\dnw\set{p}$, for all $p\in P$.
We put $\two=\set{0,1}$, the two-element poset. For quasi-ordered
sets $\seq{P,\utr_P}$ and $\seq{Q,\utr_Q}$, a map $f\colon P\to Q$
is an \emph{embedding}, if $x\utr_Py$ if{f} $f(x)\utr_Qf(y)$, for
all $x$, $y\in P$; we say that~$f$ is a \emph{lower embedding}, if
$f$ is an embedding and the range of~$f$ is an ideal of~$Q$.

For a meet-semilattice $S$, we put $S^==S\setminus\set{1}$ if $S$
has a unit, $S^==S$ otherwise. Furthermore, we denote by $\M(S)$
the set of all \mirr\ elements of $S$, that is, those $u\in S^=$
such that $u=x\wedge y$ implies that either $u=x$ or $u=y$, for
all $x$, $y\in S$. Dually, for a \jzs\ $S$, we denote by $\J(S)$
the set of all \jirr\ elements of $S$.

We denote by $\go$ the set of all natural numbers and by
$\Pow(X)$ the power set of $X$, for any set $X$.

\section{Functorial retracts}\label{S:FuncRet}

\begin{definition}\label{D:RetrAB}
Let $\cA$ and $\cB$ be
subcategories of a category $\cC$. We denote by $\Retr(\cA,\cB)$
the category whose objects and morphisms are the following:
\begin{itemize}
\item[---] \emph{Objects}: all quadruples $\seq{A,B,\eps,\mu}$,
where
$A\in\Ob\cA$, $B\in\Ob\cB$, $\eps\colon A\to\nobreak B$,
$\mu\colon B\to A$, and $\mu\circ\eps=\id_A$.

\item[---] \emph{Morphisms}: a morphism from $\seq{A,B,\eps,\mu}$
to $\seq{A',B',\eps',\mu'}$ is a pair $\seq{f,g}$, where
$f\colon A\to A'$ in $\cA$, $g\colon B\to B'$ in $\cB$,
$g\circ\eps=\eps'\circ f$, and $\mu'\circ g=f\circ\mu$ (see
Figure~\ref{Fig:RetrAB}). Composition of morphisms is defined by
the rule
$\seq{f',g'}\circ\seq{f,g}=\seq{f'\circ f,g'\circ g}$.

\end{itemize} In short, $\Retr(\cA,\cB)$ is the category of all
retractions of an object of $\cB$ onto an object of $\cA$.

The \emph{projection functor} from $\Retr(\cA,\cB)$ to $\cA$ is
the functor from $\Retr(\cA,\cB)$ to~$\cA$ that sends any object
$\seq{A,B,\eps,\mu}$ to $A$ and any morphism
$\seq{f,g}$ to $f$.
\end{definition}

\begin{figure}[htb]
 \[
 {
 \def\labelstyle{\displaystyle}
 \xymatrix{
 B\ar[r]^g\ar@{->>}[d]<.5ex>^{\mu} &
 B'\ar@{->>}[d]<.5ex>^{\mu'}\\
 A\ar@<.5ex>[u]^{\eps}\ar[r]_f & A'\ar@<.5ex>[u]^{\eps'}
 }
 }
 \]
\caption{Morphisms in $\Retr(\cA,\cB)$.}
\label{Fig:RetrAB}
\end{figure}

\begin{definition}\label{D:FunctRetr}
We say that $\cA$ is a
\emph{functorial retract} of $\cB$, if the projection functor from
$\Retr(\cA,\cB)$ to~$\cA$ has a right inverse. We shall call such
an inverse a \emph{functorial retraction} of $\cA$ to $\cB$.
\end{definition}

Hence a functorial retraction may be viewed as a triple
$\seq{\Phi,\eps,\mu}$ that satisfies the following conditions:
\begin{itemize}
\item[---] $\Phi$ is a functor from $\cA$ to $\cB$.

\item[---] For every morphism $f\colon X\to Y$ in $\cA$, we have
$\eps_X\colon X\to\Phi(X)$, $\mu_X\colon\Phi(X)\to\nobreak X$,
$\mu_X\circ\eps_X=\id_X$, $\Phi(f)\circ\eps_X=\eps_Y\circ f$, and
$\mu_Y\circ\Phi(f)=f\circ\mu_X$ (see Figure~\ref{Fig:FunctRetr}).
\end{itemize}
Observe that we do not require the diagram of
Figure~\ref{Fig:FunctRetr} to be commutative, for example,
$\Phi(f)\neq\eps_Y\circ f\circ\mu_X$ in general.

\begin{figure}[htb]
 \[
 {
 \def\labelstyle{\displaystyle}
 \xymatrix{
 \Phi(X)\ar[r]^{\Phi(f)}\ar@{->>}[d]<.5ex>^{\mu_X} &
 \Phi(Y)\ar@{->>}[d]<.5ex>^{\mu_Y}\\
 X\ar@<.5ex>[u]^{\eps_X}\ar[r]_f & Y\ar@<.5ex>[u]^{\eps_Y}
 }
 }
 \]
\caption{Functorial retraction of $\cA$ to $\cB$.}
\label{Fig:FunctRetr}
\end{figure}

\section{Sheltering between full subcategories}\label{S:Shelter}

\begin{definition}\label{D:IdMono}
An \emph{ideal of monics} of a category $\cC$ is a subcategory
$\cM$ of $\cC$ satisfying the following conditions:
\begin{enumerate}
\item Every identity of $\cC$ belongs to $\cM$.

\item $g\circ f\in\cM$ implies that $f\in\cM$, for all morphisms
$f$ and $g$ of $\cC$ such that $g\circ f$ is defined.

\item Every morphism in $\cM$ is a monic.
\end{enumerate}
\end{definition}

Of course, the monics of $\cC$ form the largest ideal of monics of
$\cC$, while the sections of $\cC$ form the smallest ideal of
monics of $\cC$. An example of often used ideal of monics distinct
from both the class of all monics and the class of all sections is 
constructed within the category of all commutative monoids, as the
ideal of all one-to-one monoid homomorphisms $f$ that satisfy
$f(x)\leq f(y)$ $\Rightarrow$ $x\leq y$, where $x\leq y$ is an
abbreviation for $(\exists z)(x+z=y)$.

\begin{definition}\label{D:shelter}
Let $\cA$ and $\cB$ be full
subcategories of a category $\cC$ and let $\cM$ be an ideal of
monics of $\cC$. A \emph{shelter of
$\cC$ by $\cB$ with respect to $\seq{\cA,\cM}$} consists of the
following data (illustrated on Figure~\ref{Fig:Shelter1}):
\begin{enumerate}
\item A functor $\xB$ from $\cC^{\iso}$ to $\cB^{\iso}$.

\item A natural transformation $S\mapsto\eta_S$ from the identity
functor on $\cC^{\iso}$ to the functor $\xB$, such that
$\eta_S\in\cM$, for every $S\in\Ob\cC$.

\item A map that with every morphism $g\colon S\to A$, where
$S\in\Ob\cC$ and $A\in\Ob\cA\cup\Ob\cB$, associates a morphism
$g^{\xB}\colon\xB(S)\to A$ such that $g=g^{\xB}\circ\eta_S$.
\end{enumerate}

Furthermore, we require the following conditions to be satisfied:
\begin{itemize}
\item[(1)] For every isomorphism $f\colon S\to T$ in $\cC$ and
every $g\colon T\to A$, with $A\in\Ob\cA\cup\Ob\cB$,
$g^{\xB}\circ\xB(f)=(g\circ f)^{\xB}$ (see
Figure~\ref{Fig:Shelter2}(1)).

\item[(2)] For every $h\colon S\to A$ and every isomorphism
$u\colon A\to A'$ with either $A$, $A'\in\Ob\cA$ or $A$,
$A'\in\Ob\cB$, $(u\circ h)^{\xB}=u\circ h^{\xB}$ (see
Figure~\ref{Fig:Shelter2}(2)).

\end{itemize}

\begin{figure}[htb]
 \[
 {
 \def\labelstyle{\displaystyle}
 \xymatrix{
 \xB(S) & & \xB(S)\ar[r]^{\xB(f)}_{\cong} & \xB(T) & &
 \xB(S)\ar[rd]^{g^{\xB}} & \\
 S\ar[u]_{\eta_S\in\cM} & & S\ar[u]^{\eta_S}\ar[r]^{\cong}_f &
 T\ar[u]_{\eta_T} & & S\ar[u]^{\eta_S}\ar[r]_g & A
 }
 }
 \]
\caption{Data describing a shelter.}
\label{Fig:Shelter1}
\end{figure}

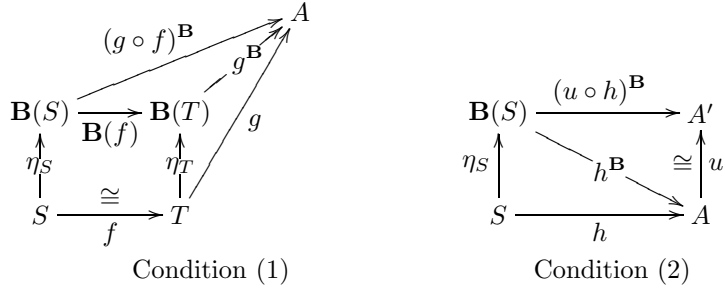
\begin{figure}[htb]
 \[
 {
 \def\labelstyle{\displaystyle}
 \xymatrix{
 & & A & & & & & \\
 \xB(S)\ar[rru]^{(g\circ f)^{\xB}}\ar[r]_{\xB(f)} &
 \xB(T)\ar[ru]|-{\ g^{\xB}} & & &
 \xB(S)\ar[rrd]|-{h^{\xB}}
 \ar[rr]^{(u\circ h)^{\xB}} & & A' \\
 S\ar[u]|-{\eta_S}\ar[r]_f^{\cong}
 \save+<15ex,-5ex>\drop{\text{Condition (1)}}\restore &
 T\ar[u]|-{\eta_T}\ar[ruu]_g & & &
 S\ar[u]^{\eta_S}\ar[rr]_h
 \save+<10ex,-5ex>\drop{\text{Condition (2)}}\restore & &
 A\ar[u]_u^{\cong}
 }
 }
 \]
\caption{Additional features of a shelter.}
\label{Fig:Shelter2}
\end{figure}

\end{definition}

\begin{remark}\label{Rk:BinA}
In all examples considered in this paper, $\cB$ is contained in
$\cA$. One can then say that a shelter is a weak reflection of
$\cC$ to $\cA$ which is everywhere a monic (i.e., all
arrows $\eta_S$ are monics), has values in $\cB$ (in case
$\cB\subseteq\cA$), and is functorial on isomorphisms.
\end{remark}

\section{Statement of the main theorem}\label{S:MainThm}

\begin{definition}\label{D:Msubobj}
Let $\cM$ be an ideal of monics of a category $\cC$. For
$S\in\Ob\cC$, we denote by $\cM(S)$ the set of all morphisms
$u\colon X\to S$ in~$\cM$, and we put\linebreak
$\cM^*(S)=\cM(S)\setminus\cM^{\iso}$. Furthermore, for
$u\colon X\to S$ and
$v\colon Y\to S$ in $\cM$, we put
 \begin{align}
 u\utr_Sv & \Longleftrightarrow
 (\exists f\colon X\to Y)(u=v\circ f);\label{Eq:ulhdv}\\
 u\sim_Sv & \Longleftrightarrow(u\utr_Sv\text{ and }v\utr_Su);
 \label{Eq:usimv}\\
 u\tr_S v & \Longleftrightarrow(u\utr_S v\text{ and }v\nutr_Su).
 \label{Eq:utrv}
 \end{align}
Obviously, $\utr_S$ is a quasi-ordering on $\cM(S)$
and $\sim_S$ is the associated equivalence. In case $u\utr_Sv$, we
shall denote by $u/v$ the unique $f\colon X\to Y$ satisfying
$u=v\circ f$. Necessarily, $f\in\cM$, and $f$ is an isomorphism
if{f} $u\sim_Sv$. We shall denote by $\lh S$ (the \emph{length of
$S$}) the length of the quasi-ordered set $\seq{\cM(S),\utr_S}$ in
case $\cM(S)$ has finite length. The blocks of $\sim_S$ will be
called the \emph{$\cM$-subobjects} of~$S$.
\end{definition}

\begin{lemma}\label{L:IncrLh}
Let $\cM$ be an ideal of monics of a category $\cC$ and let
$f\colon X\to Y$ in~$\cM$. Then the map
$\cM(f)\colon\cM(X)\to\cM(Y)$, $u\mapsto f\circ u$ is a lower
embedding. Furthermore, if both $X$ and $Y$ have finite length,
then $f$ is an isomorphism if{f} $\lh X=\lh Y$.
\end{lemma}

\begin{proof}
Verifying that $\cM(f)$ is a lower embedding is a straightforward
exercise. If $f$ is an isomorphism, then so is $\cM(f)$, thus
$\lh X=\lh Y$. If $f$ is not an isomorphism, then
$\lh X=\mathrm{height}_{\cM(X)}(\id_X)=\mathrm{height}_{\cM(Y)}(f)
<\mathrm{height}_{\cM(Y)}(\id_Y)=\lh Y$.
\end{proof}

Now we state the main technical result of the paper.

\begin{theorem}\label{T:Main}
Let $\cA$ and $\cB$ be full
subcategories of a category $\cC$ and let $\cM$ be an ideal of
monics of $\cC$. We assume the following:
\begin{enumerate}
\item Every diagram of $\cC$, indexed by a finite poset, and with
vertices either in~$\cA$ or in~$\cB$, has a colimit.

\item Every object of $\cA$ has only finitely many
$(\cA\cap\cM)$-subobjects.

\item $\cC$ is sheltered by $\cB$ with respect to $\seq{\cA,\cM}$.

\item For every $\varphi\colon A_0\to A_1$ in $\cM$ and every
\emph{section} $\eps_0\colon A_0\to B_0$, with $A_0$,
$A_1\in\Ob\cA$ and $B_0\in\Ob\cB$, there is $S\in\Ob\cC$,
together with $\eps_1\colon A_1\to S$ and
$\psi\colon B_0\to\nobreak S$ both in $\cM$, such that
$\psi\circ\eps_0=\eps_1\circ\varphi$.
\end{enumerate}

Then $\cA\cap\cM$ is a functorial retract of $\cB\cap\cM$.
\end{theorem}

{}From now on until the end of Section~\ref{S:PhiMinM}, we shall
assume that $\cA$, $\cB$, $\cC$, and~$\cM$ satisfy the assumptions
of Theorem~\ref{T:Main}, with a shelter $\xB$ denoted as in
Section~\ref{S:Shelter}. The functorial retraction of
Theorem~\ref{T:Main} will be constructed explicitly, in terms of
categorical operations and $\xB$.

Our next lemma states that in item (iv) of Theorem~\ref{T:Main},
we may assume that $S\in\Ob\cB$.

\begin{lemma}\label{L:AmalgB}
For all $\varphi\colon A_0\to A_1$ in $\cM$ and every section
$\eps_0\colon A_0\to B_0$, with $A_0$, $A_1\in\Ob\cA$ and
$B_0\in\Ob\cB$, there is $B\in\Ob\cB$, 
together with $\eps_1\colon A_1\to B$
and $\psi\colon B_0\to B$ both in~$\cM$, such that
$\psi\circ\eps_0=\eps_1\circ\varphi$.
\end{lemma}

\begin{proof}
Consider $S$, $\eps_1$, and $\psi$ obtained from
(iv) of Theorem~\ref{T:Main}. Replace $S$ by $B=\xB(S)$, $\eps_1$
by $\eta_S\circ\eps_1$, and $\psi$ by $\eta_S\circ\psi$.
\end{proof}

\begin{remark}\label{Rk:AmalgB}
Our formulation of Theorem~\ref{T:Main} is a
compromise between conciseness and generality. As one can
never be sure about future applications, let us mention a few
possible weakenings of its assumptions.
Assumption~(ii) can be weakened, by putting a cardinal upper
bound, say, $\kappa$, on the number of subobjects of all objects
of~$\cA$. Then, Assumption~(i) needs to be extended to diagrams
indexed by posets of size below~$\kappa$. However, in order to be
able to define the (ordinal) length, and, in particular, to get an
analogue of Lemma~\ref{L:IncrLh}, we need to keep the assumption
that each poset $\seq{\cM(S),\utr_S}$ is \emph{well-founded} (the
terminology \emph{artinian} is also used), that is,
every nonempty subset has a minimal element. Finally, the
diagrams involved in Assumption~(i) are fairly special, for
example, they have at least one vertex in $\cA$ and all their
arrows in~$\cM$. However, we know no situation where such
generalizations would be of any practical use.
\end{remark}

\section{Inductive construction of $\Phi$, $\eps$, $\mu$}
\label{S:ConstrPhi}

Denote by $\cA_n$ the full subcategory of $\cA$ whose objects are
those $X\in\Ob\cA$ such that $\lh X<n$, for every natural
number $n$. Of course, $\cA_0$ is the empty category.

Fix a natural number $n$, and suppose having constructed a
functor $\Phi$ from $\cA_n\cap\cM$ to $\cB$, together with a
system of morphisms $\eps_X\colon X\to\Phi(X)$ and
$\mu_X\colon\Phi(X)\to X$, for $X\in\Ob\cA_n$, such that the
following induction hypothesis is satisfied:
 \begin{equation}\label{Eq:IHPhi}
 \begin{aligned}
 \mu_X\circ\eps_X=\id_X,\ \Phi(f)\circ\eps_X=\eps_Y\circ f,
 \text{ and }\mu_Y\circ\Phi(f)=f\circ\mu_X,\\
 \text{ for every morphism }
 f\colon X\to Y\text{ in }\cA_n\cap\cM.
 \end{aligned}
 \end{equation}
\emph{We do not assume, for the moment, that $\Phi$ sends
morphisms in $\cM$ to morphisms in $\cM$}. So, for $f\colon X\to Y$
in $\cA_n\cap\cM$, all we know is that
$\Phi(f)\colon\Phi(X)\to\Phi(Y)$ in~$\cB$. We fix an object $A$ of
$\cA$ such that $\lh A=n$.

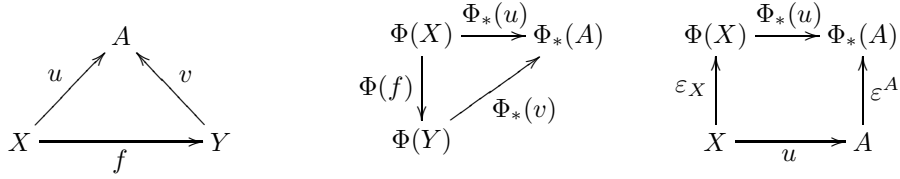
\begin{figure}[htb]
 \[
 {
 \def\labelstyle{\displaystyle}
 \xymatrix{
 & A & & & \Phi(X)\ar[r]^{\Phi_*(u)}\ar[d]_{\Phi(f)} &
 \Phi_*(A) & \Phi(X)\ar[r]^{\Phi_*(u)} & \Phi_*(A)\\
 X\ar[ru]^u\ar[rr]_f & & Y\ar[lu]_v & &
 \Phi(Y)\ar[ru]_{\Phi_*(v)} & & X\ar[u]^{\eps_X}\ar[r]_u &
 A\ar[u]_{\eps^A}
 }
 }
 \]
\caption{The colimit of $\rho_A$.}
\label{Fig:colimitrhoA}
\end{figure}

{\em Let us outline the construction. We shall introduce a
diagram $\rho_A$, indexed by a quasi-ordered set
$\seq{\AMt{A},\utr_A}$. Intuitively, $\rho_A$ consists of all
spans $\seq{u,\eta_X}$, where $X\in\Ob\cA$, $u\colon X\to A$ in
$\cM$, and $\eta_X$ is either $\eps_X$, in case
$u\notin\cM^{\iso}$, or $\id_X$, in case $u\in\cM^{\iso}$. We
equip these objects with the obvious arrows, see
\textup{\eqref{Eq:DefArrrhoA}}. An important
auxiliary construction is the colimit of $\rho_A$, which consists
of an object $\Phi_*(A)$, together with arrows
$\eps^A\colon A\to\Phi_*(A)$ and
$\Phi_*(u)\colon\Phi(X)\to\Phi_*(A)$, for $u\colon X\to A$ in
$\cM\setminus\cM^{\iso}$, subjected to the commutation relations
illustrated on Figure~\textup{\ref{Fig:colimitrhoA}}. The
resulting natural transformation $A\mapsto\eps^A$ from the
identity to $\Phi_*$ is split by $\mu^{-}\colon\Phi_*\to\id$,
$A\mapsto\mu^A$, living in
$\cC$ and given as follows: $\mu^A$ is induced by the cocone
$\seq{\id_A,\seqm{u\circ\mu_X}{u\colon X\to A
\text{ in }\cM\setminus\cM^{\iso}}}$, see
Lemma~\textup{\ref{L:A<Phi(A)}}. We observe that the definition of
$\Phi_*(A)$ does not use only all previous values of $\Phi_*$, but
really all previous values of $\Phi$.

The shelter $\xB$ is used in order to define $\Phi(A)$: namely,
$\Phi(A)=\xB(\Phi_*(A))$, see~\textup{\eqref{Eq:Phi(A)}}. The
natural transformation $A\mapsto\eps_A$, its section
$A\mapsto\mu_A$, and the arrows $\Phi(u)$, for $u\colon X\to A$ in
$\cM\setminus\cM^{\iso}$, are then defined in the natural way, see
\textup{\eqref{Eq:eps(A)}}--\textup{\eqref{Eq:Phi(u)}}. The rest
of the section is then devoted to proving that this extension
of~$\Phi$ on~$\cA_{n+1}$ can, indeed, be further extended to a
functor. Although it will turn out that $\Phi$ preserves~$\cM$,
this is a nontrivial fact and it will not be assumed as an
induction hypothesis through the construction. We shall establish
this fact in Sections~\textup{\ref{S:FactB}}
and~\textup{\ref{S:PhiMinM}}. }

Now let us go to the details.
We put $\AM{A}=(\cA\cap\cM)(A)$ and
$\AMs{A}=(\AM{A})\setminus\cM^{\iso}$. For
$u\colon X\to A$ and $v\colon Y\to A$ in $\AM{A}$ with $u\utr_Av$,
we shall often identify the morphism (in $\AM{A}$) $u\to v$ with
$f=u/v$, which is a morphism (in $\cA\cap\cM$) from~$X$ to~$Y$.
We endow the set
 \[
 \AMt{A}=\setm{\seq{u,i}\in(\AM{A})\times\two}
 {i=1\Rightarrow u\notin\cM^{\iso}}
 \]
with the partial quasi-ordering, that we shall still denote by
$\utr_A$, defined by
 \[
 \seq{u,i}\utr_A\seq{v,j}\Longleftrightarrow(u\utr_Av
 \text{ and }i\leq j),\qquad\text{for all }
 \seq{u,i},\,\seq{v,j}\in\AMt{A}.
 \]
For $\seq{u,i}\in\Ob(\AMt{A})$, where $u\colon X\to A$, we define
 \[
 \rho_A(\seq{u,i})=\begin{cases}
 X, & \text{if }i=0,\\
 \Phi(X), & \text{if }i=1.
 \end{cases}
 \]
(Observe that $i=1$ implies that $u\notin\cM^{\iso}$, thus
$\lh X<\lh A$, and thus $\Phi(X)$ is defined.) Furthermore, if
$\seq{u,i}\utr_A\seq{v,j}$ in $\AMt{A}$, we put
 \begin{equation}\label{Eq:DefArrrhoA}
 \rho_A(\seq{u,i}\to\seq{v,j})=\begin{cases}
 u/v, & \text{if }i=j=0,\\
 \Phi(u/v)\circ\eps_X, & \text{if }i=0\text{ and }j=1,\\
 \Phi(u/v), & \text{if }i=j=1.
 \end{cases}
 \end{equation}

\begin{lemma}\label{L:rhoAfunct}
The correspondence $\rho_A$ defines a functor from $\AMt{A}$ to
$\cM$.
\end{lemma}

\begin{proof}
It is obvious that $\rho_A(\seq{u,i}\to\seq{v,j})$
is a morphism from $\rho_A(\seq{u,i})$ to $\rho_A(\seq{v,j})$, and
that $\rho_A$ sends identities to identities. Now let
$\seq{u,i}\utr_A\seq{v,j}\utr_A\seq{w,k}$ in $\AMt{A}$, we need to
verify the equality
 \begin{equation}\label{Eq:rhoAassoc}
 \rho_A(\seq{u,i}\to\seq{w,k})=
 \rho_A(\seq{v,j}\to\seq{w,k})\circ\rho_A(\seq{u,i}\to\seq{v,j}).
 \end{equation}
Let $u\colon X\to A$, $v\colon Y\to A$, and
$w\colon Z\to A$, put $f=u/v$ and $g=v/w$. We separate cases.

\case{1} $i=j=k=0$. Then
 \[
 \rho_A(\seq{v,j}\to\seq{w,k})\circ\rho_A(\seq{u,i}\to\seq{v,j})
 =g\circ f=\rho_A(\seq{u,i}\to\seq{w,k}).
 \]

\case{2} $i=j=0$, $k=1$. Then
 \begin{align*}
 \rho_A(\seq{v,j}\to\seq{w,k})\circ\rho_A(\seq{u,i}\to\seq{v,j})
 &=\Phi(g)\circ\eps_Y\circ f\\
 &=\Phi(g)\circ\Phi(f)\circ\eps_X\\
 &=\Phi(g\circ f)\circ\eps_X\\
 &=\rho_A(\seq{u,i}\to\seq{w,k}).
 \end{align*}

\case{3} $i=0$, $j=k=1$. Then
 \begin{align*}
 \rho_A(\seq{v,j}\to\seq{w,k})\circ\rho_A(\seq{u,i}\to\seq{v,j})
 &=\Phi(g)\circ\Phi(f)\circ\eps_X\\
 &=\Phi(g\circ f)\circ\eps_X\\
 &=\rho_A(\seq{u,i}\to\seq{w,k}).
 \end{align*}

\case{4} $i=j=k=1$. Then
 \[
 \rho_A(\seq{v,j}\to\seq{w,k})\circ\rho_A(\seq{u,i}\to\seq{v,j})
 =\Phi(g)\circ\Phi(f)=\Phi(g\circ f)=\rho_A(\seq{u,i}\to\seq{w,k}).
 \]
This concludes the proof.
\end{proof}

\begin{lemma}\label{L:colimrhoA}
The functor $\rho_A$ has a colimit in $\cC$.
\end{lemma}

\begin{proof}
It follows from Assumption (ii) of Theorem~\ref{T:Main} that
$\AM{A}$ is equivalent to a finite poset; hence $\AMt{A}$ is also
equivalent to a finite poset. Since the colimit is a categorical
concept, the conclusion follows from Assumption~(i) of
Theorem~\ref{T:Main}.
\end{proof}

A colimit of $\rho_A$ is given by an object $\Phi_*(A)$, together
with a system of morphisms
$\theta_{\seq{u,i}}\colon\rho_A(\seq{u,i})\to\Phi_*(A)$, for all
$\seq{u,i}\in\AMt{A}$, subjected to certain commutation relations.
In case $u\in\AMs{A}$, the equality
$\theta_{\seq{u,0}}=\theta_{\seq{u,1}}\circ\eps_{\dom u}$ holds.
Hence, putting $\Phi_*(u)=\theta_{\seq{u,1}}$ and
$\eps^A=\theta_{\seq{\id_A,0}}$, we obtain that the colimit of
$\rho_A$ is given by the object $\Phi_*(A)$, together with
morphisms $\Phi_*(u)\colon\Phi(\dom u)\to\Phi_*(A)$, for
$u\in\AMs{A}$, and $\eps^A\colon A\to\Phi_*(A)$, subjected to the
commutativity of the diagrams represented on
Figure~\ref{Fig:colimitrhoA} and the universality of $\Phi_*(A)$
together with the system of morphisms consisting of all
$\Phi_*(u)$-s and $\eps^A$. Observe that for $n=0$, this reduces
to the universality of $\eps^A\colon A\to\Phi_*(A)$; so, in that
case, we may take $\Phi_*(A)=A$ and $\eps^A=\id_A$.

\begin{lemma}\label{L:A<Phi(A)}
There exists a unique morphism
$\mu^A\colon\Phi_*(A)\to A$ such that
$\mu^A\circ\nobreak\Phi_*(u)=u\circ\mu_{\dom u}$, for every
$u\in\AMs{A}$, and $\mu^A\circ\eps^A=\id_A$.
\end{lemma}

\begin{proof}
We put $\tau_u=u\circ\mu_{\dom u}$, for all
$u\in\AMs{A}$. By the universality of the colimit, it suffices to
verify that the diagrams of Figure~\ref{Fig:RetrSAA} commute, for
all $u\colon X\to A$ and $v\colon Y\to A$ in $\AMs{A}$ and $f=u/v$.
\begin{figure}[htb]
 \[
 {
 \def\labelstyle{\displaystyle}
 \xymatrix{
 \Phi(X)\ar[r]^{\tau_u}\ar[d]_{\Phi(f)} &
 A & \Phi(X)\ar[r]^{\tau_u} & A\\
 \Phi(Y)\ar[ru]_{\tau_v} & & X\ar[u]^{\eps_X}\ar[ru]_u
 }
 }
 \]
\caption{Retracting $\Phi_*(A)$ onto $A$.}
\label{Fig:RetrSAA}
\end{figure}
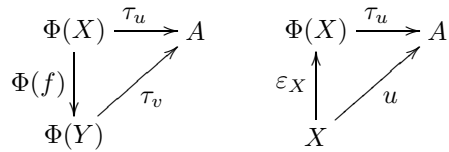

Left hand side diagram:
$\tau_v\circ\Phi(f)=v\circ\mu_Y\circ\Phi(f)= v\circ
f\circ\mu_X=u\circ\mu_X=\tau_u$.

Right hand side diagram:
$\tau_u\circ\eps_X=u\circ\mu_X\circ\eps_X=u$. This concludes the
proof.
\end{proof}

Now we are ready to define $\Phi(A)$, $\eps_A$, $\mu_A$, and
$\Phi(u)$, for $u\colon X\to A$ in $\AMs{A}$:
 \begin{align}
 \Phi(A)&=\xB(\Phi_*(A));\label{Eq:Phi(A)}\\
 \eps_A&=\eta_{\Phi_*(A)}\circ\eps^A;\label{Eq:eps(A)}\\
 \mu_A&=(\mu^A)^{\xB}.\label{Eq:mu(A)}\\
 \Phi(u)&=\eta_{\Phi_*(A)}\circ\Phi_*(u).\label{Eq:Phi(u)}
 \end{align}
These maps are represented on Figure~\ref{Fig:Phi(A)}.
They satisfy the relations
$\mu^A\circ\eps^A=\mu_A\circ\eps_A=\id_A$,
$\eps_A=\eta_{\Phi_*(A)}\circ\eps^A$,
$\mu^A=\mu_A\circ\eta_{\Phi_*(A)}$, and
$\Phi(u)=\eta_{\Phi_*(A)}\circ\Phi_*(u)$.

\begin{figure}[htb]
 \[
 {
 \def\labelstyle{\displaystyle}
 \xymatrix{
 & & \Phi(A)=\xB(\Phi_*(A))
 \ar@{->>}[dll]<-1ex>_(.6){\mu_A=(\mu^A)^{\xB}}\\
 A\ar[urr]_(.4){\eps_A}\ar@<-.5ex>[rr]_{\eps^A} & &
 \Phi_*(A)\ar[u]|-{\eta_{\Phi_*(A)}}
 \ar@{->>}[ll]<-.3ex>_(.3){\mu^A}\\
 X\ar[u]^u\ar@<-.5ex>[rr]_{\eps_X} & &
 \Phi(X)\ar[u]|-{\Phi_*(u)}\ar@{->>}[ll]<-.3ex>_(.3){\mu_X}
 \ar@/_3pc/[uu]_{\Phi(u)}
 }
 }
 \]
\caption{Defining $\Phi(A)$, $\eps_A$, $\mu_A$, and $\Phi(u)$.}
\label{Fig:Phi(A)}
\end{figure}

\begin{figure}[htb]
 \[
 {
 \def\labelstyle{\displaystyle}
 \xymatrix{
 & A & & \Phi(X)\ar[rr]^{\Phi_*(g\circ u)}\ar[d]_{\Phi(f)} & &
 \Phi_*(A') & \Phi(X)\ar[rr]^{\Phi_*(g\circ u)} & & \Phi_*(A')\\
 X\ar[ru]^u\ar[rr]_f & & Y\ar[lu]_v &
 \Phi(Y)\ar[rru]_{\Phi_*(g\circ v)} & & &
 X\ar[u]^{\eps_X}\ar[rr]_u & &
 A\ar[u]^{\eps^{A'}\circ g}
 }
 }
 \]
\caption{Putting $\Phi_*(A')$ above the diagram defining
$\Phi_*(A)$.}
\label{Fig:colimitrhoAA'}
\end{figure}
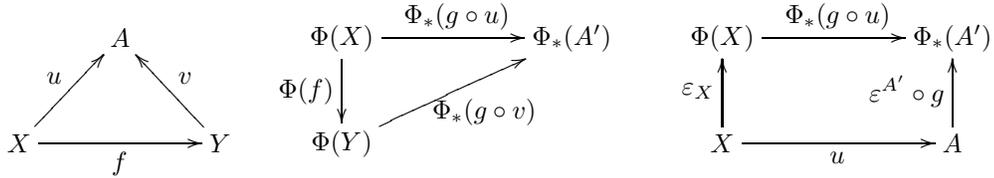

The computations of the relations
$\Phi(u)\circ\eps_X=\eps_A\circ u$ and
$\mu_A\circ\Phi(u)=u\circ\mu_X$ can be followed on
Figures~\ref{Fig:colimitrhoA} and~\ref{Fig:Phi(A)}:
 \begin{gather*}
 \Phi(u)\circ\eps_X=\eta_{\Phi_*(A)}\circ\Phi_*(u)\circ\eps_X
 =\eta_{\Phi_*(A)}\circ\eps^A\circ u=\eps_A\circ u,\\
 \mu_A\circ\Phi(u)=\mu_A\circ\eta_{\Phi_*(A)}\circ\Phi_*(u)
 =\mu^A\circ\Phi_*(u)=u\circ\mu_X.
 \end{gather*} In order to complete the extension of $\Phi$ to all
morphisms, it remains to define $\Phi(g)$, where $g\colon A\to A'$
in $\cA\cap\cM$ and $\lh A=\lh A'=n$. Observe that, by
Lemma~\ref{L:IncrLh}, $g$ is an isomorphism. Moreover, if
$u\colon X\to A$ and $v\colon Y\to A$ belong to $\AMs{A}$ with
$u\utr_Av$ and putting $f=u/v$, the diagrams of
Figure~\ref{Fig:colimitrhoAA'} commute.

Therefore, by the universal property of $\Phi_*(A)$ and the
associated limiting morphisms, there exists a unique morphism
$\ol{g}\colon\Phi_*(A)\to\Phi_*(A')$ such that
 \begin{equation}\label{Eq:ol(g)}
 \ol{g}\circ\eps^A=\eps^{A'}\circ g\text{ and }
 \Phi_*(g\circ u)=\ol{g}\circ\Phi_*(u),\text{ for all }
 u\in\AMs{A}.
 \end{equation}
Symmetrically, there exists a unique morphism
$\ol{g}'\colon\Phi_*(A')\to\Phi_*(A)$ such that
 \begin{equation}\label{Eq:ol(g')}
 \ol{g}'\circ\eps^{A'}=\eps^A\circ g^{-1}\text{ and }
 \Phi_*(u)=\ol{g}'\circ\Phi_*(g\circ u),\text{ for all }
 u\in\AMs{A}.
 \end{equation}
Again by using the universal property defining
$\Phi_*(A)$ and
$\Phi_*(A')$, we obtain that~$\ol{g}$ and~$\ol{g}'$ are mutually
inverse isomorphisms. We define
 \begin{equation}\label{Eq:Phi(g)}
 \Phi_*(g)=\ol{g}\text{ and }\Phi(g)=\xB(\Phi_*(g)).
 \end{equation}
So $\Phi_*(g)$ is an isomorphism from $\Phi_*(A)$ onto
$\Phi_*(A')$, and, since $\xB$ is a functor from~$\cC^{\iso}$ to
$\cB^{\iso}$, $\Phi(g)$ is an isomorphism from $\Phi(A)$ onto
$\Phi(A')$.

\begin{lemma}\label{L:muAolg}
$\mu^{A'}\circ\Phi_*(g)=g\circ\mu^A$.
\end{lemma}

\begin{proof}
By using the universal property defining
$\Phi_*(A)$, it suffices to verify that the diagram represented on
Figure~\ref{Fig:diagk} commutes, in both cases
$h=\mu^{A'}\circ\Phi_*(g)$ and
$h=g\circ\mu^A$, for all $u\colon X\to A$ in $\AMs{A}$. Of course,
none of the arrows of Figure~\ref{Fig:diagk} except $\eps^A$, $g$,
and $h$ are needed in case $n=0$, in which case $\Phi_*(A)=A$ and
$\eps^A=\id_A$.
\begin{figure}[htb]
 \[
 {
 \def\labelstyle{\displaystyle}
 \xymatrix{
 & & & A'\\
 \Phi(X)\ar[rr]|-{\Phi_*(u)}
 \ar@/^1pc/[rrru]|-{\mu^{A'}\circ\Phi_*(g\circ u)}
 & & \Phi_*(A)\ar[ru]|-{\ h} &\\
 X\ar[u]^{\eps_X}\ar[rr]_u & &
 A\ar[u]^{\eps^A}\ar@/_1pc/[ruu]_g &
 }
 }
 \]
\caption{Characterizing a morphism from $\Phi_*(A)$ to $A'$.}
\label{Fig:diagk}
\end{figure}
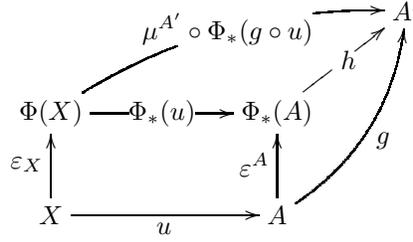
The details of the computations use \eqref{Eq:ol(g)},
\eqref{Eq:Phi(g)}, and Lemma~\ref{L:A<Phi(A)}; they are as follows:
\begin{align*}
\mu^{A'}\circ\Phi_*(g)\circ\Phi_*(u)& =\mu^{A'}\circ\Phi_*(g\circ
u),\\ g\circ\mu^A\circ\Phi_*(u)&=g\circ u\circ\mu_X
=\mu^{A'}\circ\Phi_*(g\circ u),\\
\mu^{A'}\circ\Phi_*(g)\circ\eps^A& =\mu^{A'}\circ\eps^{A'}\circ
g=g,\\ g\circ\mu^A\circ\eps^A&=g.
\end{align*}
This completes the proof.
\end{proof}

\begin{lemma}\label{L:Phi(g)circepsA}
$\Phi(g)\circ\eps_A=\eps_{A'}\circ g$ and
$\mu_{A'}\circ\Phi(g)=g\circ\mu_A$.
\end{lemma}

\begin{proof}
By using \eqref{Eq:eps(A)}, \eqref{Eq:ol(g)}, and
\eqref{Eq:Phi(g)}, we obtain
 \[
 \Phi(g)\circ\eps_A=\Phi(g)\circ\eta_{\Phi_*(A)}\circ\eps^A
 =\eta_{\Phi_*(A')}\circ\Phi_*(g)\circ\eps^A
 =\eta_{\Phi_*(A')}\circ\eps^{A'}\circ g=\eps_{A'}\circ g.
 \]
Furthermore, since $\xB$ is a shelter, we obtain, by using
\eqref{Eq:mu(A)}, the following equalities:
 \begin{align*}
 \mu_{A'}\circ\Phi(g)&=(\mu^{A'})^{\xB}\circ\xB(\Phi_*(g))
 =(\mu^{A'}\circ\Phi_*(g))^{\xB},\\
 g\circ\mu_A&=g\circ(\mu^A)^{\xB}=(g\circ\mu^A)^{\xB}.
 \end{align*}
Therefore, by Lemma~\ref{L:muAolg},
$\mu_{A'}\circ\Phi(g)=g\circ\mu_A$.
\end{proof}

At this stage, we have extended $\Phi$ to
$\cA_{n+1}\cap\cM$, up to verification of preservation of
composition by $\Phi$. Proving this preservation is the object of
the next three lemmas.

\begin{lemma}\label{L:Phin+1n+1}
Let $A_0$, $A_1$, $A_2\in\Ob\cA$ with
$\lh A_0=\lh A_1=\lh A_2=n$, let $f\colon A_0\to\nobreak A_1$ and
$g\colon A_1\to A_2$ in $\cM$. Then
$\Phi_*(g\circ f)=\Phi_*(g)\circ\Phi_*(f)$ and
$\Phi(g\circ f)=\Phi(g)\circ\Phi(f)$.
\end{lemma}

\begin{proof}
By Lemma~\ref{L:IncrLh}, both $f$ and $g$ are
isomorphisms. By using \eqref{Eq:ol(g)} and \eqref{Eq:Phi(g)}, we
obtain
 \[
 \Phi_*(g)\circ\Phi_*(f)\circ\eps^{A_0}
 =\Phi_*(g)\circ\eps^{A_1}\circ f=\eps^{A_2}\circ g\circ f,
 \]
and, for all $u\in\AMs{A_0}$,
 \[
 \Phi_*(g\circ f\circ u)=\Phi_*(g)\circ\Phi_*(f\circ u)
 =\Phi_*(g)\circ\Phi_*(f)\circ\Phi_*(u).
 \]
Since these properties determine $\Phi_*(g\circ f)$, we obtain
that $\Phi_*(g\circ f)=\Phi_*(g)\circ\Phi_*(f)$. Since $\xB$ is a
functor from $\cC^{\iso}$ to $\cB^{\iso}$, the equality
$\Phi(g\circ f)=\Phi(g)\circ\Phi(f)$ follows.
\end{proof}

\begin{lemma}\label{L:Phinn+1}
Let $X\in\Ob\cA_n$ and let $A$, $A'\in\Ob\cA$ such that
$\lh A=\lh A'=n$, let $u\colon X\to A$ and let $g\colon A\to A'$
in $\cM$. Then $\Phi(g\circ u)=\Phi(g)\circ\Phi(u)$.
\end{lemma}

\begin{proof}
By Lemma~\ref{L:IncrLh}, $g$ is an isomorphism. By using
\eqref{Eq:Phi(u)}, \eqref{Eq:ol(g)}, and \eqref{Eq:Phi(g)}, we
obtain
 \begin{multline*}
 \Phi(g\circ u)=\eta_{\Phi_*(A')}\circ\Phi_*(g\circ u)
 =\eta_{\Phi_*(A')}\circ\Phi_*(g)\circ\Phi_*(u)\\
 =\Phi(g)\circ\eta_{\Phi_*(A)}\circ\Phi_*(u)=\Phi(g)\circ\Phi(u),
 \end{multline*} which concludes the proof.
\end{proof}

\begin{lemma}\label{L:Phinnn+1}
Let $X$, $Y\in\Ob\cA_n$ and let $A\in\Ob\cA$ such that
$\lh A=n$, let $f\colon X\to\nobreak Y$ and let $u\colon Y\to A$ in
$\cM$. Then $\Phi(u\circ f)=\Phi(u)\circ\Phi(f)$.
\end{lemma}

\begin{proof}
By using \eqref{Eq:Phi(u)} and the relations on
Figure~\ref{Fig:colimitrhoA}, we obtain
 \[
 \Phi(u)\circ\Phi(f)=\eta_{\Phi_*(A)}\circ\Phi_*(u)\circ\Phi(f)
 =\eta_{\Phi_*(A)}\circ\Phi_*(u\circ f)=\Phi(u\circ f),
 \]
which concludes the proof.
\end{proof}

At this stage, $\Phi$, $\eps$, and $\mu$ have been extended to the
whole category $\cA_{n+1}\cap\cM$. Therefore, arguing by induction
on $n$, we obtain an extension of $\Phi$, $\eps$, and $\mu$ on
$\cA\cap\cM$ that satisfies the following:
 \begin{equation}\label{Eq:fullextPhi}
 \begin{aligned}
 \mu_X\circ\eps_X=\id_X,\ \Phi(f)\circ\eps_X=\eps_Y\circ f,
 \text{ and }\mu_Y\circ\Phi(f)=f\circ\mu_X,\\
 \text{ for every morphism }
 f\colon X\to Y\text{ in }\cA\cap\cM.
 \end{aligned}
 \end{equation}

\begin{definition}\label{D:CanBcov}
The triple $\seq{\Phi,\eps,\mu}$ thus constructed is the
\emph{canonical $\cB$-cover} of~$\cA$.
\end{definition}

The construction $\seq{\Phi,\eps,\mu}$ involves the shelter $\xB$
and categorical operations such as the colimit. Hence, even for
fixed $\xB$, it is defined uniquely only if we choose
representatives for colimits of diagrams: otherwise, it is defined
only up to isomorphism.

What is still missing is that we do not know yet whether the image
under $\Phi$ of a morphism in $\cM$ is a morphism in $\cM$ (which
is why we have, so far, kept this condition out of the induction
hypothesis). This is the hardest part of the proof, and it will be
the object of the next two sections.

\section{Factoring $\xB$-liftings}\label{S:FactB}

In this section we shall establish (see Lemma~\ref{L:UnivPhi}) a
certain ``quasi-universality'' property of the canonical
$\cB$-cover $\seq{\Phi,\eps,\mu}$ of $\cA$, with respect to the
notion of\linebreak \emph{$\cB$-lifting} introduced in the
following definition.

\begin{definition}\label{D:Blift}
Let $A\in\Ob\cA$, let $\cI$ be an ideal of $\AM{A}$, and let
$\Psi\colon\cI\to\cB$ be a functor. A \emph{$\cB$-lifting} of
$\Psi$ is a natural transformation from the domain functor
$u\mapsto\dom u$ (from $\cI$ to $\cA\cap\cM$) to $\Psi$.
\end{definition}

\begin{figure}[htb]
 \[
 {
 \def\labelstyle{\displaystyle}
 \xymatrix{
 & A & & \Psi(u)\ar[rr]^{\Psi(f)} & & \Psi(v)\\
 X\ar[ur]^u\ar[rr]^f & & Y\ar[ul]_v &
 X\ar[u]^{\epst_u}\ar[rr]^f & & Y\ar[u]_{\epst_v}
 }
 }
 \]
\caption{Illustrating a $\cB$-lifting of $\Psi$.}
\label{Fig:Blift}
\end{figure}
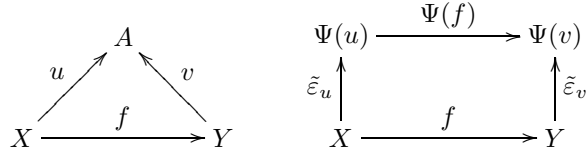

Hence a $\cB$-lifting of $\Psi$ consists of a family
$\epst=\seqm{\epst_u}{u\in\cI}$, where
$\epst_u\colon\dom u\to\nobreak\Psi(u)$, for all $u\in\cI$, such
that if $u\colon X\to A$ and $v\colon Y\to A$ in $\cI$ with
$u\utr_Av$, then, putting $f=u/v$, the equality
$\Psi(f)\circ\epst_u=\epst_v\circ f$ holds, see
Figure~\ref{Fig:Blift}. Observe that we use the convention,
introduced at the beginning of Section~\ref{S:ConstrPhi}, to
identify $f$ with $u\to v$, so $\Psi(f)$ is, in fact, defined as
$\Psi(u\to v)$.

\begin{definition}\label{D:FactorLift}
Let $A\in\Ob\cA$, let $\cI$ be an ideal of $\AM{A}$, and let
$\Psi\colon\cI\to\cB$ be a functor. A \emph{factor} of
$\seq{\Psi,\epst}$ is a natural transformation
$\delta=\seqm{\delta_u}{u\in\cI}$ from $\Phi\circ\dom$ to $\Psi$
such that $\epst_u=\delta_u\circ\eps_{\dom u}$, for all $u\in\cI$.
\end{definition}

Hence, for $u\colon X\to A$ and $v\colon Y\to A$ in $\cI$ with
$u\utr_Av$ and putting $f=u/v$, the diagram of
Figure~\ref{Fig:FactorPsi} commutes.
\begin{figure}[htb]
 \[
 {
 \def\labelstyle{\displaystyle}
 \xymatrix{
 \Psi(u)\ar[rr]^{\Psi(f)} & & \Psi(v)\\
 \Phi(X)\ar[u]_{\delta_u}\ar[rr]^{\Phi(f)} & &
 \Phi(Y)\ar[u]^{\delta_v}\\
 X\ar[u]_{\eps_X}\ar[rr]^f\ar@/^2pc/[uu]^{\epst_u} & &
 Y\ar[u]^{\eps_Y}\ar@/_2pc/[uu]_{\epst_v}
 }
 }
 \]
\caption{Illustrating a factor of $\seq{\Psi,\epst}$.}
\label{Fig:FactorPsi}
\end{figure}

\begin{lemma}\label{L:UnivPhi}
Let $A\in\Ob\cA$ and let $\cI$ and
$\cJ$ be ideals of $\AM{A}$ such that $\cJ$ contains~$\cI$. Let
$\epst$ be a $\cB$-lifting of a functor $\Psi\colon\cJ\to\cB$.
Then any factor of
$\seq{\Psi\res_{\cI},\epst\res_{\cI}}$ can be extended to a factor
of $\seq{\Psi,\epst}$.
\end{lemma}

\begin{proof}
Arguing by induction on the length reduces the
problem to the case where $\setm{v\in\AM{A}}{v\tr_Au}$ is
contained in $\cI$, for all
$u\in\cJ$. So let $u\colon U\to A$ in $\cJ\setminus\cI$, we shall
define a morphism $\delta_u\colon\Phi(U)\to\Psi(u)$.

For all $v\colon X\to U$ and $w\colon Y\to U$ in $\AMs{U}$
such that $v\utr_Uw$, letting $f=v/w$, both
relations $u\circ v\tr_Au$ and $u\circ w\tr_Au$ hold, thus both
$u\circ v$ and $u\circ w$ belong to $\cI$. Furthermore, the
diagrams of Figure~\ref{Fig:getdelu} commute: this is obvious for
the left hand side, while for the right hand side,
$\Psi(v)\circ\delta_{u\circ v}\circ\eps_X
=\Psi(v)\circ\epst_{u\circ v}=\epst_u\circ v$.

\begin{figure}[htb]
 \[
 {
 \def\labelstyle{\displaystyle}
 \xymatrix{
 \Phi(X)\ar[d]_{\Phi(f)}\ar[rr]^{\delta_{u\circ v}} & &
 \Psi(u\circ v)\ar[d]_{\Psi(f)}\ar[r]^{\Psi(v)} & \Psi(u) &
 \Phi(X)\ar[rr]^{\Psi(v)\circ\delta_{u\circ v}} & & \Psi(u)\\
 \Phi(Y)\ar[rr]^{\delta_{u\circ w}} & &
 \Psi(u\circ w)\ar[ru]_{\Psi(w)} & &
 X\ar[u]_{\eps_X}\ar[rr]^v & & U\ar[u]_{\epst_u}
 }
 }
 \]
\caption{Putting $\Psi(u)$ above the diagram defining $\Phi_*(U)$.}
\label{Fig:getdelu}
\end{figure}

Hence, by the universal property defining $\Phi_*(U)$, there
exists a unique morphism $\gamma_u\colon\Phi_*(U)\to\Psi(u)$ such
that $\gamma_u\circ\Phi_*(v)=\Psi(v)\circ\delta_{u\circ v}$, for
all $v\in\AMs{U}$, and $\epst_u=\gamma_u\circ\eps^U$. Put
$\delta_u=\gamma_u^{\xB}$ (see Figure~\ref{Fig:Defdeltau}).
\begin{figure}[htb]
 \[
 {
 \def\labelstyle{\displaystyle}
 \xymatrix{
 \Phi_*(U)\ar[r]^{\gamma_u} & \Psi(u) &
 \Phi_*(U)\ar[r]^{\gamma_u} & \Psi(u) &
 \Phi(U)\ar[rd]^(.6){\delta_u=(\gamma_u)^{\xB}} &\\
 \Phi(\dom v)\ar[u]_{\Phi_*(v)}\ar[r]_{\delta_{u\circ v}} &
 \Psi(u\circ v)\ar[u]_{\Psi(v)} &
 U\ar[u]^{\eps^U}\ar[ru]_{\epst_u} & &
 \Phi_*(U)\ar[u]^{\eta_{\Phi_*(U)}}\ar[r]_{\gamma_u} & \Psi(u)
 }
 }
 \]
\caption{Defining $\gamma_u$ and $\delta_u$.}
\label{Fig:Defdeltau}
\end{figure}

We verify that the $\delta_u$-s are as required (see
Figure~\ref{Fig:FactorPsi}). First,
 \[
 \delta_u\circ\eps_U=\delta_u\circ\eta_{\Phi_*(U)}\circ\eps^U
 =\gamma_u\circ\eps^U=\epst_u.
 \]
Our next series of calculations will prove that the extended
$\delta$ is a natural transformation from $\Phi\circ\dom$
to~$\Psi$.

For $u\colon U\to A$ in $\cJ\setminus\cI$ and $v\colon X\to U$ in
$\AMs{U}$ (thus in $\cI$),
 \[
 \delta_u\circ\Phi(v)=\delta_u\circ\eta_{\Phi_*(U)}\circ\Phi_*(v)
 =\gamma_u\circ\Phi_*(v)=\Psi(v)\circ\delta_{u\circ v}.
 \]
Now let $u\colon U\to A$ and $v\colon V\to A$ in
$\cJ\setminus\cI$ such that $u\utr_Av$, and
put $f=u/v$. If $f$ is not an isomorphism, then
$u\tr_Av$, thus (since $v\in\cJ$) $u\in\cI$, a contradiction. Hence
$f$ is an isomorphism. We prove that
$\Psi(f^{-1})\circ\gamma_v\circ\Phi_*(f)$ satisfies the properties
defining $\gamma_u$.
 \begin{align*}
 \Psi(f^{-1})\circ\gamma_v\circ\Phi_*(f)\circ\eps^U
 &=\Psi(f^{-1})\circ\gamma_v\circ\eps^V\circ f &&
 (\text{by }\eqref{Eq:ol(g)}\text{ and }\eqref{Eq:Phi(g)})\\
 &=\Psi(f^{-1})\circ\epst_v\circ f\\
 &=\epst_u.
 \end{align*}
Let $w\in\AMs{U}$. By using \eqref{Eq:ol(g)} and \eqref{Eq:Phi(g)},
we compute:
 \begin{align*}
 \Psi(f^{-1})\circ\gamma_v\circ\Phi_*(f)\circ\Phi_*(w)
 &=\Psi(f^{-1})\circ\gamma_v\circ\Phi_*(f\circ w)\\
 &=\Psi(f^{-1})\circ\Psi(f\circ w)\circ\delta_{u\circ w}\\
 &=\Psi(w)\circ\delta_{u\circ w}.
 \end{align*}
Therefore,
$\Psi(f^{-1})\circ\gamma_v\circ\Phi_*(f)=\gamma_u$, that is,
$\gamma_v\circ\Phi_*(f)=\Psi(f)\circ\gamma_u$. Now we can compute
further, using the assumption that $\xB$ is a shelter:
 \begin{align*}
 \delta_u=\gamma_u^{\xB}
 &=\bigl(\Psi(f^{-1})\circ\gamma_v\circ\Phi_*(f)\bigr)^{\xB}\\
 &=\Psi(f^{-1})\circ\gamma_v^{\xB}\circ\xB(\Phi_*(f))\\
 &=\Psi(f^{-1})\circ\delta_v\circ\Phi(f),
 \end{align*}
whence $\Psi(f)\circ\delta_u=\delta_v\circ\Phi(f)$.
\end{proof}

\section{Preservation of $\cM$ by $\Phi$}\label{S:PhiMinM}

In this section we shall prove the remaining claim about the
canonical $\cB$-cover $\seq{\Phi,\eps,\mu}$, namely, that
$\Phi$ preserves~$\cM$. The idea of the proof is the following.
For $f\colon U\to A$ in $\cA\cap\cM$, we construct, using
the amalgamation property stated in Lemma~\ref{L:AmalgB}, a certain
functor $\Psi$, defined on all subobjects of~$A$, together with a
$\cB$-lifting~$\epst$ of~$\Psi$. Furthermore, we shall see that
the restriction of $\seq{\Psi,\epst}$ to all subobjects of~$A$
below~$f$ has a factor. By the ``quasi-universality'' property
established in Section~\ref{S:FactB}, namely,
Lemma~\ref{L:UnivPhi}, this factor extends to a factor of
$\seq{\Psi,\epst}$ on all subobjects of~$A$. As $\Psi$ is
constructed in such a way that the arrow $\Psi(f\to\id_A)$ belongs
to~$\cM$, it follows that $\Phi(f)$ also belongs to~$\cM$.

\begin{lemma}\label{L:PhitoEmb}
For any morphism $f$ of $\cA\cap\cM$, $\Phi(f)$ belongs to
$\cB\cap\cM$.
\end{lemma}

\begin{proof}
We let $f\colon U\to A$ in $\cA\cap\cM$, we prove that
$\Phi(f)\in\cM$. If $f$ is an isomorphism, then (since $\Phi$ is a
functor) so is $\Phi(f)$, thus $\Phi(f)\in\cM$. {}From now on we
assume that~$f$ is not an isomorphism.

Since $\eps_U\colon U\to\Phi(U)$ is a section (for
$\mu_U\circ\eps_U=\id_U$), it follows from Lemma~\ref{L:AmalgB}
that there exists $B\in\Ob\cB$, together with
$\varphi\colon\Phi(U)\to B$ and $\eps\colon A\to B$ in $\cM$, such
that $\varphi\circ\eps_U=\eps\circ f$.

For each $u\colon X\to A$ in $\AM{A}$, we define
$\Psi(u)\in\Ob\cB$ and $\epst_u\colon X\to\Psi(u)$ by
 \begin{equation}\label{Eq:DefPsiepst}
 \Psi(u)=\begin{cases}
 \Phi(X),&\text{if }u\utr_Af,\\
 B,&\text{otherwise};
 \end{cases}\qquad
 \epst_u=\begin{cases}
 \eps_X,&\text{if }u\utr_Af,\\
 \eps\circ u,&\text{otherwise}.
 \end{cases}
 \end{equation}
For $u\colon X\to A$ and $v\colon Y\to A$ in
$\AM{A}$, we put $g=u/v$ and we define a morphism
$\Psi(u\to v)\colon\Psi(u)\to\Psi(v)$ in $\cB$ as follows:

\case{1} $v\utr_Af$. Put $\Psi(u\to v)=\Phi(g)$.

\case{2} $u\utr_Af$, $v\nutr_Af$. Put
$\Psi(u\to v)=\varphi\circ\Phi(u/f)$.

\case{3} $u,v\nutr_Af$. Put $\Psi(u\to v)=\id_B$.

\setcounter{claim}{0}
\begin{claim}\label{Cl:PsiEpst}
In the context above,
$\Psi(u\to v)\circ\epst_u=\epst_v\circ g$.
\end{claim}

\begin{cproof}
In Case~1, this is equivalent to the statement
$\Phi(g)\circ\eps_X=\eps_Y\circ g$, which holds.

In Case~2, putting $\ol{u}=u/f$, we compute
 \[
 \varphi\circ\Phi(\ol{u})\circ\eps_X=\varphi\circ\eps_U\circ\ol{u}
 =\eps\circ f\circ\ol{u}=\eps\circ u=\eps\circ v\circ g,
 \]
which is the desired statement.

In Case~3, from $v\circ g=u$ it follows that
$\eps\circ v\circ g=\eps\circ u$, which is the desired statement.
\end{cproof}

\begin{claim}\label{Cl:PsiFunct}
$\Psi$ is a functor.
\end{claim}

\begin{cproof}
It suffices to prove that
$\Psi(u\to w)=\Psi(v\to w)\circ\Psi(u\to v)$, for all
$u\utr_Av\utr_Aw$ in $\AM{A}$. Put $X=\dom u$, $Y=\dom v$,
$Z=\dom w$, $g=u/v$, and $h=v/w$. We separate cases.

\case{1} $w\utr_Af$. Then
 \[
 \Psi(v\to w)\circ\Psi(u\to v)=\Phi(h)\circ\Phi(g)=
 \Phi(h\circ g)=\Psi(u\to w).
 \]

\case{2} $v\utr_Af$ and $w\nutr_Af$. Put $\ol{u}=u/f$ and
$\ol{v}=v/f$. The equality $u=v\circ g$ can be written
$f\circ\ol{u}=f\circ\ol{v}\circ g$, thus, since $f$ is monic,
$\ol{u}=\ol{v}\circ g$. Therefore,
 \[
 \Psi(u\to w)=\varphi\circ\Phi(\ol{u})=
 \varphi\circ\Phi(\ol{v})\circ\Phi(g)=
 \Psi(v\to w)\circ\Psi(u\to v).
 \]

\case{3} $u\utr_Af$ and $v\nutr_Af$. Put $\ol{u}=u/f$. Then
 \[
 \Psi(u\to w)=\varphi\circ\Phi(\ol{u})=
 \id_B\circ\varphi\circ\Phi(\ol{u})=
 \Psi(v\to w)\circ\Psi(u\to v).
 \]

\case{4} $u\nutr_Af$. Then
$\Psi(u\to v)=\Psi(u\to w)=\Psi(v\to w)=\id_B$, whence
 \[
 \Psi(u\to w)=\Psi(v\to w)\circ\Psi(u\to v).
 \]
This concludes the proof of our claim.
\end{cproof}

By Claims~\ref{Cl:PsiEpst} and \ref{Cl:PsiFunct},
$\epst$ is a $\cB$-lifting of the functor $\Psi$ on the ideal
$\cJ=\AM{A}$. Furthermore, putting
 \[
 \cI=\setm{u\in\AM{A}}{u\utr_Af},
 \]
the rule $\delta_u=\id_{\Phi(\dom u)}$ defines a factor of
$\seq{\Psi\res_{\cI},\epst_{\cI}}$. Therefore, by
Lemma~\ref{L:UnivPhi}, $\delta$ extends to a factor of
$\seq{\Psi,\epst}$, which we shall still denote by $\delta$.

Since $f\in\cI$ and $\id_A\notin\cI$, we obtain
$\Psi(f\to\id_A)=\varphi\circ\Phi(\id_U)=\varphi$. Therefore,
 \[
 \varphi=\varphi\circ\delta_f=\delta_{\id_A}\circ\Phi(f).
 \]
Since $\varphi\in\cM$, we obtain that $\Phi(f)\in\cM$.
\end{proof}

This completes the proof of Theorem~\ref{T:Main}: the canonical
$\cB$-cover $\seq{\Phi,\eps,\mu}$ is a solution of the given
problem.

\section{Distributive and ultraboolean semilattices}
\label{S:DistrSem}

A \jzs\ $\seq{S,\vee,0,\leq}$ is \emph{distributive}, if its ideal
lattice $\Id S$ is distributive, see G. Gr\"atzer
\cite[Section~II.5]{GLT2}. Equivalently, for all $a$, $b$, $c\in
S$, if $c\leq a\vee b$, then there are $x\leq a$ and $y\leq b$
such that $c=x\vee y$. Distributive \jzs s are characterized in P.
Pudl\'ak \cite[Fact~4, p.~100]{Pudl85} as directed \jz-unions of
finite distributive \jzs s.

Denote by $\cS$ the category of finite \jzs s and \jzh s, and by
$\cD$ and $\cB$ the full subcategories of $\cS$ consisting of all
distributive, respectively Boolean members of $\cS$. We denote by
$\cM$ the subcategory of $\cS$ consisting of all \jze s. Of
course, $\cM$ is an ideal of monics of $\cS$.

\begin{lemma}\label{L:CanExt}
Let $S$, $T$, and $D$ be \jzs s with $D$ finite distributive, and
let $e\colon S\to T$ and $f\colon S\to D$ be \jzh s with $e$ an
\emph{embedding}. Then there exists a largest \jzh\
$g\colon T\to D$ extending $f$, and it is given by the formula
 \begin{equation}\label{Eq:DefExt}
 g(t)=\bigwedge\nolimits^D\bigl(f(s)\mid s\in S,\ t\leq e(s)\bigr),
 \text{ for all }t\in T.
 \end{equation}
\pup{Of course, the meet of the empty set is defined here as the
unit of $D$.}
\end{lemma}

\begin{proof}
An easy exercise. Although the distributivity of $D$ is not used
for correctness of the definition \eqref{Eq:DefExt}, it is
used for proving that $g$ is a \jzh.
\end{proof}

We shall call the map $g$ defined in \eqref{Eq:DefExt} the
\emph{largest extension} of $f$ with respect to $e$.

Any $S\in\Ob\cS$ is a finite lattice, thus also a meet-semilattice.
We put $\xB(S)=\Pow(\M(S))$, and we let
 \[
 \eta_S\colon S\into\xB(S),\ a\mapsto\setm{u\in\M(S)}{a\nleq u}.
 \]

\begin{lemma}[folklore]\label{L:EmbSBS}
The map $\eta_S$ is a \jzue\ from $S$ into
$\seq{\xB(S),\cup,\es,S}$, for every finite \jzs\ $S$.
\end{lemma}

For an isomorphism $f\colon S\to T$ of finite \jzs s, we put
 \[
 \xB(f)(X)=f[X],\text{ for all }X\in\xB(S).
 \]
It is immediate to verify that $\xB$ is a functor from
$\cS^{\iso}$ to $\cB^{\iso}$ and that $\eta$ is a natural
transformation from the identity of $\cS^{\iso}$ to $\xB$.

\begin{definition}\label{D:SheltSem}
Let $S$ and $A$ be finite \jzs s with $A$ distributive, and let
$g\colon S\to A$ be a \jzh. We denote $g^{\xB}$ the largest
\jz-extension of~$g$ from $\xB(S)$ to $A$
with respect to the embedding $\eta_S\colon S\into\xB(S)$.
\end{definition}

Now the proof of the following lemma is a straightforward exercise
(see Definition~\ref{D:shelter}). Items~(i), (ii), (1), and (2) of
Definition~\ref{D:shelter} follows from the fact that the formulas
defining $\xB$, $g\mapsto g^{\xB}$, and $\eta$ are `intrinsic',
thus preserved under isomorphisms. Item~(iii) of
Definition~\ref{D:shelter} follows from Lemma~\ref{L:CanExt} (see
Definition~\ref{D:SheltSem}).

\begin{lemma}\label{L:SemShelt}
The correspondences $\xB$, $\eta$ described above define a
shelter of $\cS$ by $\cB$ with respect to $\seq{\cD,\cM}$.
\end{lemma}

The corresponding commutative diagram is given on
the right hand side of Figure~\ref{Fig:Shelter1}. Now we are ready
to prove our main semilattice-theoretical result.

\begin{theorem}\label{T:FRSemil}
There exists a functorial retraction $\seq{\Phi,\eps,\mu}$ of the
category $\cD\cap\cM$ to the category $\cB\cap\cM$. Furthermore,
$\eps_A$ is a \jzue, for all $A\in\Ob\cD$.
\end{theorem}

\begin{proof}
We prove that the assumptions of Theorem~\ref{T:Main} are
satisfied, where we replace $\cC$ by $\cS$ and $\cA$ by $\cD$.
Item~(i) is a very particular case of the well-known fact that
every variety of algebras has small colimits, see, for
example, \cite[Theorem~8.3.8]{Berg}. Item~(ii) is trivial.
Item~(iii) is Lemma~\ref{L:SemShelt}.

Finally, it is proved in \cite[Theorem~2.10]{HoKi71} that every
semilattice embeds into an injective semilattice. Hence the
variety of semilattices has the so-called Transfer Property (see
\cite[Proposition~1.5]{KMPT}), thus \emph{a fortiori} the
Amalgamation Property. Since every finitely generated
semilattice is finite, these results extend to the finite case.
Technically speaking, these results are
established in \cite{HoKi71,KMPT} for semilattices which do not
necessarily have a unit; however, the extension to the case with
unit is trivial. The result for \jzs s is dual. This obviously
implies Assumption~(iv) of Theorem~\ref{T:Main}.

It remains to establish that the maps $\eps_A$ constructed in the
proof of Theorem~\ref{T:Main} are $1$-preserving. We argue by
induction on $\lh A$. By definition, $\Phi_*(A)$ is \jz-generated
by the set
 \begin{equation}\label{Eq:GenPhi*(A)}
 G=\im\eps^A\cup\bigcup\famm{\im\Phi_*(u)}{u\in\AMs{A}}.
 \end{equation}
Let $u\colon X\into A$ in $\AMs{A}$.
By the induction hypothesis, $\eps_X$ is $1$-preserving, thus
 \[
 \Phi_*(u)(1_{\Phi(X)})=\Phi_*(u)\circ\eps_X(1_X)=
 \eps^A\circ u(1_X)\leq\eps^A(1_A),
 \]
thus the largest element of $X$ is $\eps^A(1_A)$. Hence the
largest element of $\Phi_*(A)$ is also $\eps^A(1_A)$, that is,
$\eps^A$ is $1$-preserving. Since $\eta_S$ is $1$-preserving for
all $S$, it follows that $\eps_A=\eta_{\Phi_*(A)}\circ\eps^A$ is
also $1$-preserving.
\end{proof}

Finally, denote by $b_A$ the largest $b\in\Phi(A)$ such that
$\mu_A(b)=0$, for every finite distributive \jzs\ $A$. After
replacing $\Phi(A)$ by its interval $[b_A,1]$, and this for all
$A$, we obtain that the $\mu_A$-s may be assumed to separate zero
(i.e., $\mu_A^{-1}\set{0}=\set{0}$).

For convenience, we list here the properties satisfied by the
triple $\seq{\Phi,\eps,\mu}$ of Theorem~\ref{T:FRSemil}:
 \begin{itemize}
 \item The correspondence $\Phi$ is a functor from the category
 $\cD\cap\cM$ of all finite distributive \jzs s with \jze s to the
 category $\cB\cap\nobreak\cM$ of all finite Boolean \jzs s with
 \jze s.
 
 \item The map $\eps_A$ is a \jzue\ from $A$ into $\Phi(A)$
 and the map $\mu_A$ is a zero-separating \jzuh\ from $\Phi(A)$
 onto $A$, for every finite distributive \jzs\ $A$. Furthermore,
 $\mu_A\circ\eps_A=\id_A$.
 
 \item For every \jze\ $f\colon X\into Y$ between finite
 distributive \jzs s $X$ and $Y$, both equalities
 $\Phi(f)\circ\eps_X=\eps_Y\circ f$ and
 $\mu_Y\circ\Phi(f)=f\circ\mu_X$ hold.
 \end{itemize}
Observe that these properties imply that $\Phi(f)$ preserves the
unit whenever $f$ does.

As shown by the following result, this functorial inverse of the
functor~$\Pi$ does not arise from an adjunction.

\begin{proposition}\label{P:NoAdj}
The projection functor
$\Pi\colon\Retr(\cD\cap\cM,\cB\cap\cM)\to\cD\cap\cM$ has neither a
right nor a left adjoint.
\end{proposition}

\begin{proof}
Denote the objects of $\cR=\Retr(\cD\cap\cM,\cB\cap\cM)$ as
$p=\seq{D_p,B_p,\alpha_p,\beta_p}$. A left or right adjoint of
$\Pi$ is given by a functor $\Psi\colon\cD\to\cR$. Let $\Psi$ be
given, for any $f\colon D\to E$ in $\cD\cap\cM$, by
 \begin{equation}\label{Eq:DefPsi}
 \begin{aligned}
 \Psi(D)=\seq{\tilde{D},\Phi(D),\eps^D,\mu^D}\text{ and }
 \Psi(f)=\seq{\tilde{f},\Phi(f)},\\
 \text{ where }
 \tilde{f}\colon\tilde{D}\to\tilde{E}\text{ and }
 \Phi(f)\colon\Phi(D)\to\Phi(E).
 \end{aligned}
 \end{equation}
Suppose first that $\Psi$ is a left adjoint of $\Pi$, and denote
by $\eta$ the unit of the corresponding adjunction. So
$\eta_D\colon D\into\tilde{D}$, for every $D\in\Ob\cD$. By the
definition of an adjunction, for all $D\in\Ob\cD$,
$p\in\Ob\cR$, and $f\colon D\into D_p$, there exists a unique
$\seq{f_*,f^*}\colon\Psi(D)\to p$ such that $f=f_*\circ\eta_D$. In
particular, for $D=D_p$ and $f=\id_D$ (as seen in the Introduction,
there are always $B$, $\alpha$, and $\beta$ such that
$\seq{D,B,\alpha,\beta}\in\Ob\cR$), we obtain, since all our
semilattices are finite, that $\eta_D$ is an \emph{isomorphism}
from~$D$ onto~$\tilde{D}$.

Now let $D=\three$ (the three-element chain), $D_p=B_p=\two^2$,
$\alpha_p=\beta_p=\id_{\two^2}$, and $f\colon\three\into\two^2$
any \jze. {}From $\seq{f_*,f^*}\colon\Psi(D)\to p$ it follows that
$f_*\circ\mu^D=f^*$, an embedding. Thus $\mu^D$ is an embedding,
and so $\Phi(D)\cong\tilde{D}\cong D\cong\three$, which is not
Boolean; a contradiction.

Now suppose that $\Psi$ is a right adjoint of $\Pi$, and denote by
$\eps$ the counit of the corresponding adjunction. So
$\eps_D\colon\tilde{D}\into D$, for all $D\in\Ob\cD$.
By the definition of an
adjunction, for all $D\in\Ob\cD$, $p\in\Ob\cR$, and
$f\colon D_p\into D$, there exists a unique
$\seq{f_*,f^*}\colon p\to\Psi(D)$ such that $\eps_D\circ f_*=f$.
In particular, $f^*\colon B_p\into\Phi(D)$, and hence
$|B_p|\leq|\Phi(D)|$. However, for $D$ and $D_p$ fixed, $|B_p|$
can be taken arbitrarily large, a contradiction.
\end{proof}

The functorial retraction given by Theorem~\ref{T:FRSemil} is
given by an explicit formula. This makes it possible to give a
crude upper bound for the maximum $\varphi(n)$ of all
cardinalities of $\Phi(A)$, where $A$ is a distributive \jzs\ of
cardinality at most $n$. Of course, $\varphi(1)=1$.
As \eqref{Eq:GenPhi*(A)} gives a
generating subset of $\Phi_*(A)$, we obtain
 \[
 |\Phi_*(A)|\leq 2^{n+1+2^n\varphi(n)},
 \]
whence
 \[
 \varphi(n+1)\leq 2^{2^{n+1+2^n\varphi(n)}}.
 \]
Hence $\varphi(n)$ is, roughly speaking, majorized by a
tower of exponentials of length~$2n$, which is, of course, beyond
the reach of any implementation.

As illustrated in \cite{Tuma93}, the poset of distributive
subsemilattices of a finite distributive \jzs\ can be quite
complicated. To the contrary, the corresponding structure is much
nicer for Boolean subsemilattices. This motivates the definition
of ultraboolean introduced in Section~\ref{S:Intro}.

\begin{corollary}\label{C:FRSemil}
Every distributive \jzs\ $D$ is a \jz-retract of an
ultraboolean \jzs\ $B$. Furthermore, if $D$ has a unit, then $B$
can be taken with a unit.
\end{corollary}

\begin{proof}
By Pudl\'ak's Lemma, $D=\varinjlim_{i\in I}D_i$, for a direct
system $\seqm{D_i,f_{i,j}}{i\leq j\text{ in }I}$ of finite
distributive \jzs s and \jze s $f_{i,j}\colon D_i\to D_j$.
Furthermore, we may assume that all $D_i$-s contain as an element
the unit of $D$ in case there is any.

Now we use the functorial retraction constructed in the proof
of Theorem~\ref{T:FRSemil}.
The semilattice $B=\varinjlim_{i\in I}\Phi(D_i)$, with transition
maps $\Phi(f_{i,j})\colon\Phi(D_i)\to\Phi(D_j)$, is ultraboolean
and has a unit in case $D$ has a unit. Furthermore, the natural
transformations $\seqm{\eps_{D_i}}{i\in I}$ and
$\seqm{\mu_{D_i}}{i\in I}$ define, by direct limit, \jzh s
$\eps\colon D\to B$ and $\mu\colon B\to D$ such that
$\mu\circ\eps=\id_D$. Therefore, $D$ is a retract of $B$.
If $D$ has a unit, then $\eps(1_D)=1_B$.
\end{proof}

\section{Simultaneous lattice embeddings into finite Boolean
semilattices}\label{S:SimultEmb}

The maps $\eps_X\colon X\to\Phi(X)$ constructed in the proof of
Theorem~\ref{T:FRSemil} are \jzu-embeddings. On the other hand,
for every finite distributive \jzs\ $S$, the embedding from $S$
into $\Pow(\J(S))$ that with every $a\in S$ associates the set
$\setm{p\in\J(S)}{p\leq a}$ is always a lattice embedding, and it
has nice ``almost functorial'' properties, see
\cite[Section~1]{Ruzi}. Hence the question whether a new
functorial retraction may be constructed, with the corresponding
maps $\eps_X$-s being \emph{lattice} homomorphisms,
is natural.

In the present section we shall prove, by a counterexample, that
this is not possible. In fact we shall prove a much stronger
negative statement, see Example~\ref{Ex:UnliftLatt}.

All direct systems considered in this section will be indexed by
posets. Hence, if $\seq{I,\leq}$ is a poset, an $I$-indexed direct
system in a category $\cA$ consists of a system
$\seqm{A_i,f_{i,j}}{i\leq j\text{ in }I}$, with
$f_{i,j}\colon A_i\to A_j$ in $\cA$, for $i\leq j$ in $I$, such
that $f_{i,i}=\id_{A_i}$ and $f_{i,k}=f_{j,k}\circ f_{i,j}$, for
all $i\leq j\leq k$ in $I$.

Suppose now that all $A_i$-s are finite \jzs s, all the
$f_{i,j}$-s are \jze s, and let $i\leq j$
in $I$. For any $q\in\J(A_j)$, we define $\partial^{i,j}q$ as the
set of all \emph{minimal} $p\in A_i$ such that $q\leq f_{i,j}(p)$.
Of course, $\partial^{i,j}q$ is a subset of $\J(A_i)$.

\begin{definition}\label{D:FunctLE}
Let $\alpha=\seqm{A_i,f_{i,j}}{i\leq j\text{ in }I}$ and
$\beta=\seqm{B_i,g_{i,j}}{i\leq j\text{ in }I}$ be direct systems
of lattices, indexed by the same poset $I$.
A \emph{simultaneous lattice embedding} from $\alpha$ into $\beta$
is a system $\seqm{\eps_i}{i\in I}$ of lattice embeddings
$\eps_i\colon A_i\into B_i$ such that
$\eps_j\circ f_{i,j}=g_{i,j}\circ\eps_i$, for all $i\leq j$ in $I$.
\end{definition}

\begin{proposition}\label{P:NecessLattLift}
Let $\seqm{A_i,f_{i,j}}{i\leq j\text{ in }I}$ be a direct system,
indexed by a poset~$I$, of finite distributive \jzs s and \jze s,
that admits a simultaneous lattice embedding into a direct system
of finite Boolean semilattices. Then for all $i\leq j$ in $I$ and
all $p\in\J(A_i)$, there exists $q\leq f_{i,j}(p)$ in $\J(A_j)$
such that the following statements hold:
\begin{enumerate}
\item $\partial^{i,j}q=\set{p}$.

\item For all $k\in I$ with $i\leq k\leq j$ and all
$r\in\partial^{k,j}q$, the following implication holds:
 \[
 r\leq f_{i,k}(1_{A_i})\ \Longrightarrow r\leq f_{i,k}(p).
 \]
\end{enumerate}
\end{proposition}

\begin{proof}
We fix a simultaneous lattice
embedding as in Definition~\ref{D:FunctLE}, with all the $A_i$-s
finite distributive \jzs s and all the $B_i$-s finite Boolean.

We put $0_i=0_{A_i}$ and $1_i=1_{A_i}$, for all $i\in I$.
Furthermore, define $X_i$ as the set of all atoms of $B_i$.
Replacing $B_i$ by the interval $[\eps_i(0_i),\eps_i(1_i)]$ (which
is still Boolean), we see that there is no loss of generality in
assuming that $\eps_i$ is a $\seq{\vee,\wedge,0,1}$-embedding, for
all $i\in I$. Furthermore, we may assume that the $X_i$-s are
pairwise disjoint and that all the $f_{i,j}$-s are set-theoretical
inclusion mappings, so that $A_i\subseteq A_j$, for all $i\leq j$
in $I$.

Since $\eps_i$ is a $(\wedge,1)$-homomorphism, we can define
$\mu_i(x)$, for $x\in B_i$, as the least $a\in A_i$ such that
$x\leq\eps_i(a)$. Hence,
 \begin{equation}\label{Eq:Dualepsmu}
 x\leq\eps_i(a)\ \Longleftrightarrow\ \mu_i(x)\leq a,
 \text{ for all }
 \seq{a,x}\in A_i\times B_i.
 \end{equation}
For $\seq{\xi,\eta}\in X_i\times X_j$, let $\xi\utr\eta$ hold, if
$\eta\leq g_{i,j}(\xi)$. It is obvious that $\utr$ is a partial
ordering on $\bigcup\famm{X_i}{i\in I}$.

The following claim records a few elementary facts.

\setcounter{claim}{0}
\begin{claim}\label{Cl:Basicepsmu}\hfill
\begin{enumerate}
\item $\mu_i\circ\eps_i=\id_{A_i}$.

\item $x\leq\eps_i\circ\mu_i(x)$, for all $x\in B_i$.

\item $\mu_i(\xi)\in\J(A_i)$, for all $\xi\in X_i$.

\item $\eps_i(a)=\bigvee\famm{\xi\in X_i}{\mu_i(\xi)\leq a}$, for
all $a\in A_i$.

\item $\mu_j\circ g_{i,j}(x)\leq\mu_i(x)$, for all $x\in B_i$.

\item $\xi\utr\eta$ implies that $\mu_j(\eta)\leq\mu_i(\xi)$, for
all $\seq{\xi,\eta}\in X_i\times X_j$.
\end{enumerate}
\end{claim}

\begin{proof}
It follows from \eqref{Eq:Dualepsmu} that $\eps_i(a)\leq\eps_i(b)$
if{f} $\mu_i\eps_i(a)\leq b$, for all $a$, $b\in A_i$. Since
$\eps_i$ is an embedding, (i) follows.

Substituting $a=\mu_i(x)$ in \eqref{Eq:Dualepsmu} gives
immediately~(ii).

Put $p=\mu_i(\xi)$, for $\xi\in X_i$. Since $\xi$ is nonzero, so
is $p$. Let $p=a\vee b$, where $a$, $b\in A_i$. So,
$\xi\leq\eps_i(p)=\eps_i(a)\vee\eps_i(b)$, but $\xi$ is an atom of
$B_i$, whence either $\xi\leq\eps_i(a)$ or $\xi\leq\eps_i(b)$, and
hence, by the definition of $p$, either $p=a$ or $p=b$. Item (iii)
follows.

By \eqref{Eq:Dualepsmu}, an element $\xi$ of $X_i$ lies
below $\eps_i(a)$ if{f} $\mu_i(\xi)\leq a$. Since $B_i$ is
Boolean, (iv) follows.

{}From (ii) it follows that
$g_{i,j}(x)\leq g_{i,j}\circ\eps_i\circ\mu_i(x)
=\eps_j\circ\mu_i(x)$; item (v) follows.

For $\seq{\xi,\eta}\in X_i\times X_j$ with $\xi\utr\eta$, that is,
$\eta\leq g_{i,j}(\xi)$, we obtain, using (v), that
$\mu_j(\eta)\leq\mu_j\circ g_{i,j}(\xi)\leq\mu_i(\xi)$. Item (vi)
follows.
\end{proof}

\begin{claim}\label{Cl:Liftpartial}
The set $\partial^{i,j}\mu_j(\eta)$ is contained in
$\setm{\mu_i(\xi)}{\xi\in X_i\text{ and }\xi\utr\eta}$, for all
$i\leq j$ in $I$ and all $\eta\in X_j$.
\end{claim}

\begin{cproof}
Let $p\in\partial^{i,j}\mu_j(\eta)$. {}From
Claim~\ref{Cl:Basicepsmu}(iv) it follows that
 \[
 \eps_j(p)=\bigvee\famm{\beta\in X_j}{\mu_j(\beta)\leq p},
 \]
whence $\eta\leq\eps_j(p)$. Therefore, again by using
Claim~\ref{Cl:Basicepsmu}(iv), we obtain
 \[
 \eta\leq g_{i,j}\circ\eps_i(p)=
 g_{i,j}\left(\bigvee\famm{\xi\in X_i}{\mu_i(\xi)\leq p}\right)
 =\bigvee\famm{g_{i,j}(\xi)}
 {\xi\in X_i\text{ and }\mu_i(\xi)\leq p},
 \]
whence there exists $\xi\in X_i$ such that $\xi\utr\eta$ and
$\mu_i(\xi)\leq p$. Hence, by Claim~\ref{Cl:Basicepsmu}(vi),
$\mu_j(\eta)\leq\mu_i(\xi)\leq p$, with $\mu_i(\xi)\in A_i$ and
$p\in\partial^{i,j}\mu_j(\eta)$. Therefore, $p=\mu_i(\xi)$.
\end{cproof}

Now we can conclude the proof of
Proposition~\ref{P:NecessLattLift}. It follows from
Claim~\ref{Cl:Basicepsmu}(i,iv) that
 \[
 p=\mu_i\circ\eps_i(p)=
 \bigvee\famm{\mu_i(\xi)}{\xi\in X_i\text{ and }\mu_i(\xi)\leq p},
 \]
thus, since $p$ is \jirr, there exists $\xi\in X_i$ such that
$p=\mu_i(\xi)$. Since $g_{i,j}$ is an embedding, there exists
$\eta\in X_j$ such that $\eta\leq g_{i,j}(\xi)$ and
$\eta\nleq g_{i,j}(\neg^{B_i}\xi)$. This means that $\xi\utr\eta$
and $\xi'\nutr\eta$ for all $\xi'\in X_i\setminus\set{\xi}$.
We prove that the element
$q=\mu_j(\eta)$ is as desired.

First, by Claim~\ref{Cl:Basicepsmu}(iii), $q$ belongs to $\J(A_j)$.
Since $q\leq p$, the set $\partial^{i,j}q$ is nonempty. Let
$p'\in\partial^{i,j}q$. By Claim~\ref{Cl:Liftpartial}, there
exists $\xi'\in X_i$ such that $\xi'\utr\eta$ and
$\mu_i(\xi')=p'$. By the definition of $\eta$, we obtain that
$\xi'=\xi$, so $p'=p$. Hence, $\partial^{i,j}q=\set{p}$.

Now let $k\in I$ with $i\leq k\leq j$ and let
$r\in\partial^{k,j}q$ with $r\leq 1_i$. The latter inequality
implies that $\partial^{i,k}r$ is nonempty. Let
$p'\in\partial^{i,k}r$. Since
$r\in\partial^{k,j}q=\partial^{k,j}\mu_j(\eta)$, there exists, by
Claim~\ref{Cl:Liftpartial}, $\zeta\in X_k$ such that
$r=\mu_k(\zeta)$ and $\zeta\utr\eta$. Since $p'\in\partial^{i,k}r$,
there exists, again by Claim~\ref{Cl:Liftpartial}, $\xi'\in X_i$
such that $\xi'\utr\zeta$ and $\mu_i(\xi')=p'$. So $\xi'\utr\eta$
with $\xi'\in X_i$, whence $\xi'=\xi$, and so $p'=p$. Therefore,
by Claim~\ref{Cl:Basicepsmu}(vi),
$r=\mu_k(\zeta)\leq\mu_i(\xi')=p$.
\end{proof}

Now we obtain the promised counterexample.

\begin{example}\label{Ex:UnliftLatt}
There exists a square \pup{i.e., a diagram indexed by $\two^2$} of
finite distributive \jzs s and \jzue s that does not have any
simultaneous lattice embedding into any diagram of finite Boolean
semilattices.
\end{example}

\begin{proof}
Identify the finite poset $P$ diagrammed on the left hand side
of Figure~\ref{Fig:PSA1} with its canonical image in the
(distributive) lattice $A$ of all ideals of $P$. So,
$1=p_1\vee p_2$ and $P=\J(A)$. Put $p=p_1\wedge p_2=q_1\vee q_2$
and let $S$ denote the \jz-subsemilattice of $A$ generated by
$\set{p,p_1,p_2}$. Hence $S$ is distributive and
$\J(S)=\set{p,p_1,p_2}$.
 \begin{figure}[hbt]
 \includegraphics{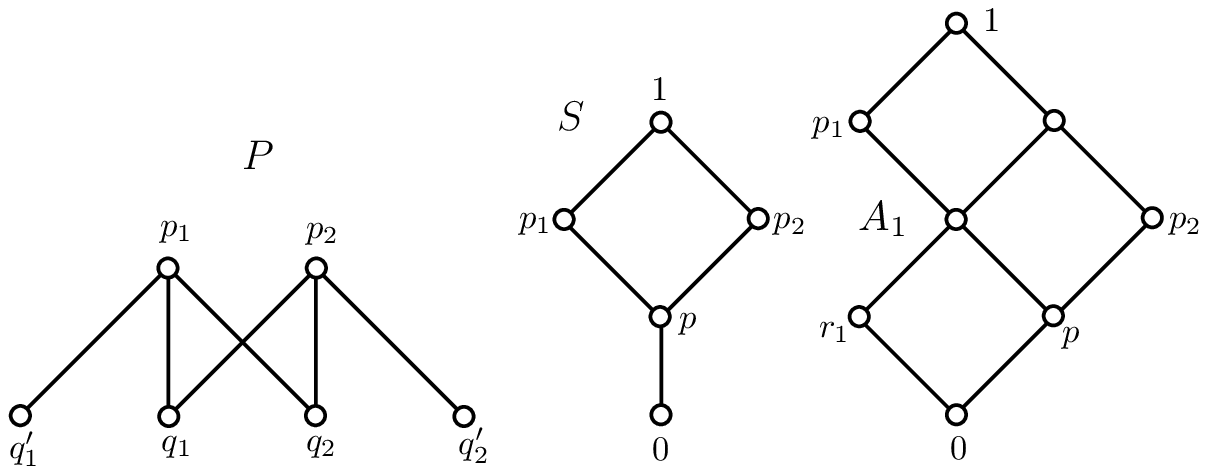}
 \caption{The poset $P$ and the \jzus s $S$ and $A_1$.}
 \label{Fig:PSA1}
 \end{figure}
For $i\in\set{1,2}$, put $r_i=q_i\vee q'_i$ and
$P_i=\set{p,p_1,p_2,r_i}$. Denote by~$A_i$ the \jz-subsemilattice
of $A$ generated by $P_i$. The only nontrivial comparable pairs
in~$P_1$ are given by $r_1<p_1$ and $p<p_1,p_2$. Furthermore, since
$p_1\nleq r_1\vee p_2$, all elements of~$P_1$
are join-prime in $A_1$, hence $A_1$ is isomorphic to the lattice
of all ideals of $P_1$; whence it is distributive. Similarly,
$A_2$ is distributive. The semilattices~$S$ and $A_1$ are
diagrammed on Figure~\ref{Fig:PSA1}. Of course, $A_1$ and $A_2$
are isomorphic.

So we have obtained four \jzus s $S\subseteq A_1,A_2\subseteq A$.
Suppose that this square has a simultaneous lattice embedding into
a diagram of finite Boolean semilattices. We apply
Proposition~\ref{P:NecessLattLift} to the element $p\in\J(S)$. The
element $q\in\nobreak\J(A)$ given by
Proposition~\ref{P:NecessLattLift} lies below $p$, so we may
assume, by symmetry, that $q=q_1$. {}From
$r_1\in\partial^{A_1,A}q_1$ and $r_1\leq 1_S$, it follows that
$r_1\leq p$, a contradiction.
\end{proof}

The proof of Proposition~\ref{P:NecessLattLift} above makes
essential use of the distributivity of all $B_i$-s. As we shall
see in Section~\ref{S:Tisch}, this is unavoidable.

\section{The Gr\"atzer-Schmidt extension and retracts
of ultra-simple-atomistic semilattices}
\label{S:Tisch}

A well-known result by M. Tischendorf \cite{Tisch} gives a direct
construction implying that \emph{every finite lattice embeds into
some finite atomistic lattice}. For a finite lattice~$L$, denote
by $\Ti(L)$ the finite atomistic lattice obtained from $L$
\emph{via} Tischendorf's construction.
It is proved in \cite{Tisch} that $\Ti(L)$ is a finite atomistic
lattice containing (as a bounded lattice) $L$, \emph{via} the
embedding $\zeta_L\colon L\into\Ti(L)$, $a\mapsto\dnw a\cap\J(L)$.
In fact, it is proved in \cite{Tisch} that $\zeta_L$ is
\emph{congruence-preserving}, that is, the natural map from
the congruence lattice $\Con L$ of $L$ to $\Con\Ti(L)$ is an
isomorphism. Furthermore, the map $\rho_L\colon\Ti(L)\to L$,
$X\mapsto\bigvee X$ is easily seen to be a \jz-retraction of
$\zeta_L$.

Although the correspondence $L\mapsto\Ti(L)$ cannot be extended
``naturally'' to arbitrary \jze s, it can be extended to
\emph{isomorphisms}. This is sufficient to construct from it an
appropriate shelter. This shelter can, in turn, be used to prove
the following analogue of Corollary~\ref{C:FRSemil}: \emph{Every
\jzs\ is a retract of a \jzs\ which is a directed \jz-union of
finite atomistic lattices}.

However, a much stronger result can be proved with a much simpler
method, see Theorem~\ref{T:GSretr}. We shall now present this
proof.

We recall that a \jzs\ $K$ is \emph{atomistic}, if every
element of $K$ is a join of atoms of $K$.
The purpose of the first part of the following definition is to
separate the two distinct notions of simple semilattice (which is
a trivial) and simple lattice.

\begin{definition}\label{D:UltraSAt}
Let $K$ be a lattice with zero. We say that the
\jzs~$\seq{K,\vee,0}$ is \emph{lattice-simple}, if the lattice
$\seq{K,\vee,\wedge}$ is simple. A \jzs~$K$ is
\emph{ultra-simple-atomistic}, if $K$ is the directed union of its
finite, lattice-simple, atomistic
\jz-subsemilattices.
\end{definition}

We denote by $\At K$ the set of all atoms of a \jzs\ $K$, and we
put $\NAt K=K\setminus(\set{0}\cup\At K)$. For every $a\in\NAt K$,
we adjoin distinct atoms $p_a^i<a$, for $i\in\set{0,1}$, such
that $p_a^i=p_b^j$ only in case $a=b$ and $i=j$. Now we put
 \[
 \GS(K)=K\cup\setm{p_a^i}{a\in\NAt K\text{ and }i<2}.
 \]
Since this construction is used in the
proof of \cite[Lemma~7]{GrSc99}, we shall call it the
\emph{Gr\"atzer-Schmidt extension} of $K$.

The ordering of $\GS(K)$ consists of the ordering of $K$,
augmented by the following pairs:
 \begin{align*}
 p_a^i\leq b&\ \Longleftrightarrow\ a\leq b
 &&(\text{in case }a\in\NAt K),\\
 a\leq p_b^j&\ \Longleftrightarrow\ a=0
 &&(\text{in case }b\in\NAt K),\\
 p_a^i\leq p_b^j&\ \Longleftrightarrow\ (a=b\text{ and }i=j)
 &&(\text{in case }a,\,b\in\NAt K),
 \end{align*}
for $a$, $b\in K$ and $i$, $j<2$.

The following lemma records a few straightforward properties of
$\GS(K)$.

\begin{lemma}\label{L:GSLatt}
Let $K$ be a \jzs. Then the following properties hold:
\begin{enumerate}
\item The ordering $\leq$ endows $\GS(K)$ with a structure of \jzs.

\item The inclusion map $\eps_K\colon K\into\GS(K)$ is a
complete $\seq{\vee,\wedge}$-embedding \pup{that is, an
order-embedding that preserves all meets and joins defined in
$K$}.

\item If $K$ is a lattice, then so is $\GS(K)$.

\item The map $\mu_K\colon\GS(K)\to K$ extending $\id_K$ such that
$\mu_K(p_a^i)=a$, for all $a\in\NAt K$ and $i<2$, is a \jzh,
and $\mu_K\circ\eps_K=\id_K$.
\end{enumerate}
\end{lemma}

For elements $a$, $b$, and $c$ in a \jzs\ $K$, we say that
$c=a\oplus b$, if $c=a\vee b$ and $a\wedge b=0$. Moreover, we say
that $a$ and $b$ are \emph{perspective}, in notation
$a\sim b$, if there exists $x\in K$ such that
$a\oplus x=b\oplus x$. The following lemma contains further related
properties.

\begin{lemma}\label{L:GSSimple}
Let $K$ be a \jzs. Then the following properties hold:
\begin{enumerate}
\item Every element of $\GS(K)$ is a join of at most two atoms.

\item For all $x$, $y\in\GS(K)$ such that $0<x<y$, there exists an
atom $p$ of $\GS(K)$ such that $y=x\oplus p$.

\item Any two atoms of $\GS(K)$ are perspective.

\item If $K$ is a lattice, then $\GS(K)$ is lattice-simple.

\item In the general case, $\GS(K)$ is ultra-simple-atomistic.
\end{enumerate}
\end{lemma}

\begin{proof}
(i) Any $a\in\NAt K$ satisfies that $a=p_a^0\vee p_a^1$.

(ii) Necessarily, $y\in K$.
If $x\in K$, then $y=x\oplus p_y^0$. If $x=p_a^i$, then
$a\leq y$, and so $p_a^i\oplus p_y^{1-i}=y$.

(iii) Let $x$ and $y$ be distinct atoms of $\GS(K)$, and put
$c=\mu_K(x)\vee\mu_K(y)$. If $x$, $y\in\At K$, then
$x\oplus p_c^0=y\oplus p_c^0=c$. If $x\in\At K$ and $y=p_a^i$
(so $\mu_K(y)=a$), then $x\oplus p_c^{1-i}=y\oplus p_c^{1-i}=c$.
Suppose that $x=p_a^i$ and $y=p_b^j$. Since $a\in\NAt K$, there
exists $d\in K$ such that $0<d<a$. If $a=b$, then
$x\oplus d=y\oplus d=a$. Suppose that $a\neq b$, say $b\nleq a$.
Then $x\oplus p^{1-j}_c=y\oplus p^{1-j}_c=c$.

(iv) It follows from Lemma~\ref{L:GSLatt}(iii) that
$\GS(K)$ is a lattice.
Denote by $\Theta(x,y)$ the (lattice-)congruence of $\GS(K)$
generated by the pair $\seq{x,y}$, for any $x$, $y\in\GS(K)$.
It follows from (iii) that $\Theta(0,x)=\Theta(0,y)$, for all
atoms $x$ and $y$ of $\GS(K)$. Therefore, by (i) (or (ii)),
$\GS(K)$ is lattice-simple.

(v) As $K$ is the directed union of its finite \jz-subsemilattices
(we define the empty directed union as $\set{0}$),
we obtain that $\GS(K)$ is the directed union of all $\GS(F)$,
for $F$ a nontrivial finite join-subsemilattice of $K$. By (i)
and (iv), $\GS(F)$ is finite, atomistic, and lattice-simple, for
all such $F$.
\end{proof}

For \jzs s $K$ and $L$ and a \jze\ $f\colon K\into L$, we define a
map $\GS(f)\colon\GS(K)\to\GS(L)$ by the rule
 \begin{align*}
 \GS(f)(a)&=a&&\text{ for }a\in K,\\
 \GS(f)(p_a^i)&=p_{f(a)}^i&&\text{ for }a\in\NAt K
 \text{ and }i<2.
 \end{align*}
The verification of the following lemma is straightforward.

\begin{lemma}\label{L:GSfunctor}
In the context above, the map $\GS(f)$ is a \jze\ from $\GS(K)$
into $\GS(L)$ such that $\GS(f)\circ\eps_K=\eps_L\circ f$ and
$\mu_L\circ\GS(f)=f\circ\mu_K$. Furthermore, if $f$ is a lattice
homomorphism, then so is $\GS(f)$.
\end{lemma}

Putting together some of the information above, we obtain the
following rather elementary result.

\begin{theorem}\label{T:GSretr}
The triple $\seq{\GS,\eps,\mu}$ is a functorial retraction of the
category of \jzs s and \jze s to the full subcategory of
ultra-simple-atomistic \jzs s. Furthermore, the functor $\GS$
sends finite lattices to finite \pup{lattice-}simple lattices.
\end{theorem}

The essence of this result can be captured by the following
somewhat loose formulation: \emph{Every \jzs\ is a retract of some
ultra-simple-atomistic \jzs, and this holds functorially}.

Further properties of the functorial retraction of
Theorem~\ref{T:GSretr} are obtained above. For example,
for any \jze\ $f\colon K\into L$,
\begin{itemize}
\item[(1)] the map $\eps_K$ is a complete
$\seq{\vee,\wedge}$-embedding (this is why the assumption of
distributivity of the $B_i$-s is unavoidable in the proof of
Proposition~\ref{P:NecessLattLift});

\item[(2)] if $f\colon K\into L$ is a $0$-lattice embedding, then
so is $\GS(f)$.
\end{itemize}

Let us keep the notation of Section~\ref{S:DistrSem} for $\cS$ and
$\cM$, and denote by $\cS_{\at}$ the full subcategory of atomistic
members of~$\cS$. We state the following analogue of
Proposition~\ref{P:NoAdj}.

\begin{proposition}\label{P:NoAdjTisch}
The projection functor from
$\Retr(\cS\cap\cM,\cS_{\at}\cap\cM)$ to $\cS\cap\cM$ has neither a
right nor a left adjoint.
\end{proposition}

The proof of Proposition~\ref{P:NoAdjTisch} is
virtually the same as the one of Proposition~\ref{P:NoAdj}. 

However, as shows the following easy result
and since there are finite non-atomistic lattices, the analogue of
Theorem~\ref{T:GSretr} for \emph{lattices} does not hold.

\begin{proposition}\label{P:NonLatt}
Any finite $\seq{\vee,\wedge}$-homomorphic
image of a $\seq{\vee,\wedge}$-direct limit of finite atomistic
lattices is atomistic.
\end{proposition}

\begin{proof}
Let $K$ be a finite lattice and let $g\colon L\onto K$ be a
surjective lattice homomorphism, where $L=\varinjlim_{i\in I}L_i$,
with $I$ directed, the lattices $L_i$ finite atomistic, and
transition maps $f_i\colon L_i\to L$. Since $K$ is finite, there
exists $i\in I$ such that $g_i=g\circ f_i$ is surjective. Since
$g_i(p)$ is an atom of $K$, for any atom $p$ of $L_i$, $K$ is
atomistic.
\end{proof}

\section{Open problems}\label{S:Pbs}

As observed above, the functorial retraction constructed in the
proof of Theorem~\ref{T:FRSemil}, although theoretically
computable, lives \emph{a priori} beyond the reach of any
implementation. The most natural question is thus whether such a
functorial retraction could be constructed with `reasonable'
growth.

A possible way to formulate this problem is the following. We use
the notation of Section~\ref{S:DistrSem}.

\begin{problem}\label{Pb:GrowthPhi}
Are there a functorial retraction $\seq{\Phi,\eps,\mu}$ of
$\cD\cap\cM$ to $\cB\cap\cM$ such that $|\J(\Phi(D))|$ is bounded
by a polynomial in $|\J(D)|$, for every finite distributive
(semi)lattice $D$?
\end{problem}

Both Example~\ref{Ex:UnliftLatt} and the huge upper bound for the
construction of Theorem~\ref{T:FRSemil} suggest that $\Phi(D)$
needs to be large with respect to $D$.

Say that a \jzh\ $\mu\colon S\to T$ is \emph{weakly distributive},
if whenever $\mu(c)=a\vee b$, there is a decomposition $c=x\vee y$
in $S$ such that $\mu(x)\leq a$ and $\mu(y)\leq b$.
In view of some lifting results with respect to
the congruence functor on lattices (see \cite{CLPSurv} for a
survey), the following problem may be relevant.

\begin{problem}\label{Pb:WDRetr}
Is every distributive \jzs\ a weakly distributive image of some
ultraboolean \jzs?
\end{problem}

We know that Problem~\ref{Pb:WDRetr} has a positive answer for
\emph{countable} semilattices.

We do not know whether the analogue of Corollary~\ref{C:FRSemil}
for \emph{dimension groups} holds. By definition, a partially
ordered abelian group $G$ is a \emph{dimension group}, if~$G$ is
directed (for its ordering), unperforated, and has the
interpolation property, see \cite{Gpoag}. Special cases of
dimension groups are the \emph{simplicial groups}, that is, those
partially ordered abelian groups that are
isomorphic to some finite power of the integers (with
componentwise ordering). As defined in \cite{CaWe}, a partially
ordered abelian group is \emph{$\mathbf{E}$-ultrasimplicial}, if it
is a directed union of simplicial groups. Every
$\mathbf{E}$-ultrasimplicial group is a dimension group; the
converse is easily seen to be false, even in the divisible case
(see \cite[Example~1.2]{CaWe}).

\begin{problem}\label{Pb:Ultrasimp}
Is every dimension group a retract of some
$\mathbf{E}$-ultrasimplicial group?
\end{problem}

A similar question can be formulated in the context of \cite{GPW}.
We denote by $\Rep$ (resp., $\Rep^*$) the class of all monoids
which are direct limits (resp., directed unions) of finite
products of monoids of the form $(\ZZ/n\ZZ)\cup\set{0}$ for
positive integers $n$. A first-order
characterization of
$\Rep$ is obtained in \cite{GPW}.

\begin{problem}\label{Pb:GPWRetr}
Is every member of $\Rep$ a retract of some member of $\Rep^*$?
\end{problem}

The result of Theorem~\ref{T:FRSemil} is made possible by the
shelter $\xB$. In order to define a shelter we need a functor
playing the role of the `functor from~$\cS^{\iso}$
to~$\cB^{\iso}$' of Definition~\ref{D:shelter}(i). A special
feature of such functors is that they can be easily defined on
\emph{isomorphisms} (because they are given by `explicit'
constructions), but not on
\emph{embeddings}. There are probably many such objects within
mathematical practice. For example, it is proved in
\cite[Theorem~1.11]{AGT}, \emph{via} an explicit construction, that
every finite \jsd\ lattice embeds into some finite atomistic \jsd\
lattice. (A lattice is \emph{\jsd}, if it satisfies the
quasi-identity
$x\vee y=x\vee z\Rightarrow x\vee y=x\vee(y\wedge z)$.) This
suggests the following problem.

\begin{problem}\label{Pb:SD+}
Say that a \jzs\ $S$ is \emph{\jsd}, if for all $a$, $b$, $c\in
S$, if $a\vee b=a\vee c$, then there exists $x\leq b,c$ such that
$a\vee b=a\vee x$. Is every \jsd\ \jzs\ a retract of
some direct limit of finite atomistic \jsd\ \jzs s?
\end{problem}

\section*{Acknowledgment}

This work was partially completed while the author was visiting
the Charles University (Prague). Excellent conditions provided by
the Department of Algebra are greatly appreciated.


\begin{thebibliography}{99}
\bibitem{AGT}
K.\,V. Adaricheva, V.\,A. Gorbunov, and V.\,I. Tumanov,
\emph{Join-semidistributive lattices and convex geometries},
Adv. Math. \textbf{173} (2003), 1--49.

\bibitem{Berg}
G.\,M. Bergman,
``An Invitation to General Algebra and Universal Constructions'',
pub. Henry Helson, 15 the Crescent, Berkeley, CA, 94708, 1998.
398~p. Available online at
\texttt{http://math.berkeley.edu/\~{}gbergman/}.

\bibitem{CaWe}
J.\,F. Caillot and F. Wehrung,
\emph{Finitely presented, coherent, and ultrasimplicial ordered
abelian groups}, Semigroup Forum \textbf{61} (2000), 116--137.

\bibitem{Gpoag}
K.\,R. Goodearl,
``Partially Ordered Abelian Groups with Interpolation'',
Math. Surveys and Monographs \textbf{20}, Amer. Math. Soc.,
Providence, 1986.

\bibitem{GPW}
K.\,R. Goodearl, E. Pardo, and F. Wehrung,
\emph{Semilattices of groups and inductive limits of Cuntz
algebras}, J. Reine Angew. Math., to appear.

\bibitem{GLT2}
G. Gr\"atzer, ``General Lattice Theory. Second
edition'', new appendices by the author with B.\,A. Davey, R.
Freese, B. Ganter, M. Greferath, P. Jipsen, H.\,A. Priestley, H.
Rose, E.\,T. Schmidt, S.\,E. Schmidt, F. Wehrung, and R. Wille.
Birkh\"auser Verlag, Basel, 1998. xx+663~p.

\bibitem{GrSc99}
G. Gr{\"a}tzer and E.\,T. Schmidt,
\emph{Congruence-preserving extensions of finite lattices to
sectionally complemented lattices}, Proc. Amer. Math. Soc.
\textbf{127} (1999), 1903--1915.

\bibitem{HoKi71}
A. Horn and N. Kimura,
\emph{The category of semilattices},
Algebra Universalis \textbf{1} (1971), 26--38.

\bibitem{KMPT}
E.\,W. Kiss, L. M\'arki, P. Pr\"ohle, and W. Tholen,
\emph{Categorical algebraic properties: a compendium on
amalgamation, congruence extension, epimorphisms, residual
smallness, and injectivity},
Studia Sci. Math. Hungar. \textbf{2} (1983), 79--141.

\bibitem{McLa}
S. Mac Lane, ``Categories for the Working
Mathematician'', 2nd ed. Graduate Texts in Mathematics \textbf{5},
Springer-Verlag, New-York, Berlin, Heidelberg, 1998. xii+314~p.

\bibitem{P5}
P.\,P. P\'alfy and P. Pudl\'ak,
\emph{Congruence lattices of finite algebras and intervals in
subgroup lattices of finite groups}, Algebra Universalis
\textbf{11} (1980), 22--27.

\bibitem{Pudl85}
P. Pudl\'ak,
\emph{On congruence lattices of lattices}, Algebra Universalis
\textbf{20} (1985), 96--114.

\bibitem{Ruzi}
P. R\r{u}\v{z}i\v{c}ka,
\emph{Lattices of two-sided ideals of locally matricial algebras
and the $\Gamma$-invariant problem}, Israel J. Math. \textbf{142}
(2004), 1--28.

\bibitem{Tisch}
M. Tischendorf,
\emph{The representation problem for algebraic distributive
lattices}, Ph.D. thesis, TH Darmstadt, 1992.

\bibitem{Tuma93}
J. T\r{u}ma,
\emph{On simultaneous representations of distributive lattices},
Acta Sci. Math. (Szeged) \textbf{58} (1993), 67--74.

\bibitem{CLPSurv}
J. T\r{u}ma and F. Wehrung,
\emph{A survey of recent results on congruence lattices of lattices},
Algebra Universalis \textbf{48}, no.~4 (2002), 439--471.

\bibitem{RetLift}
F. Wehrung,
\emph{Lifting retracted diagrams with respect to projectable
functors}, preprint.

\end{thebibliography}
\end{document}